\documentclass[11pt]{article}
\title{Types of graded tensor products of graded von Neumann algebras}
\author{Jumpei Tanaka }

\usepackage{amsmath, amssymb}
\usepackage{amsfonts}
\usepackage{amsthm}
\usepackage{mathtools}
\mathtoolsset{showonlyrefs=true}
\theoremstyle{definition}
\newtheorem{dfn}{Definition}[section]

\newtheorem{prop}[dfn]{Proposition}
\newtheorem{lem}[dfn]{Lemma}
\usepackage{mathrsfs}

\begin{document}
\maketitle
\tableofcontents
\begin{abstract}
There is a famous multiplication table of types of tensor product of two von Neumann algebras. We filled out the multiplication table of graded tensor product of two graded von Neumann algebras in special cases.
\end{abstract}
\section{introduction}
It is famous fact that, the types of tensor product of two von Neumann algebras are as in the multiplication table below. In this paper, we consider the types of graded tensor product of graded\ von\ Neumann\ algebras.
In section 2, we check the definitions of some kinds of graded von Neumann algebras, and introduce some theorems we use in this paper.\ In section $3.1$, we prove that if at least one graded von Neumann algebra is of type I\hspace{-1.2pt}I\hspace{-1.2pt}I, the types of graded tensor product is equal to that of normal cases. In section $3.2$, we consider the case both $(\mathcal{R} _1,\ \text{Ad}_{\Gamma_1})$,\ $(\mathcal{R}_2,\
\text{Ad}_{\Gamma_2})$ are central and both $\mathcal{R}_1$,\ $\mathcal{R}_2$ are not
factors. Consequently, we proved they are equal to types of normal tensor\ product except the case both is of type I. In the case both is of type\ I, 
the index is twice as large as the index of normal tensor\ product. In section $3.3$, we study the case $(\mathcal{R}_2,\ \text{Ad}_{\Gamma_2})$ is central, $\mathcal{R}_1$ is a
factor, and $\mathcal{R}_2$ is not a factor. We get multiplication table which is the same as that of normal tensor product. In section $3.4$, we introduce other facts.

\begin{table}[]
    \begin{tabular}{cccccc}
                                               &
\multicolumn{5}{c}{type of $\mathcal{S}$}

                  \\ \cline{2-6}
    \multicolumn{1}{c|}{type of $\mathcal{R}$} &
\multicolumn{1}{c|}{I$_n$, n:finite} & \multicolumn{1}{c|}{I$_n$,
n:infinite} & \multicolumn{1}{c|}{I\hspace{-1.2pt}I$_1$}        &
\multicolumn{1}{c|}{I\hspace{-1.2pt}I$_{\infty}$} & I\hspace{-1.2pt}I\hspace{-1.2pt}I \\ \hline
    \multicolumn{1}{c|}{I$_m$, m:finite}       &
\multicolumn{1}{c|}{I$_{mn}$}        & \multicolumn{1}{c|}{I$_{mn}$}
       & \multicolumn{1}{c|}{I\hspace{-1.2pt}I$_1$}        &
\multicolumn{1}{c|}{I\hspace{-1.2pt}I$_{\infty}$} & I\hspace{-1.2pt}I\hspace{-1.2pt}I \\ \hline
    \multicolumn{1}{c|}{I$_m$, m:infinite}     &
\multicolumn{1}{c|}{I$_{mn}$}        & \multicolumn{1}{c|}{I$_{mn}$}
       & \multicolumn{1}{c|}{I\hspace{-1.2pt}I$_{\infty}$} &
\multicolumn{1}{c|}{I\hspace{-1.2pt}I$_{\infty}$} & I\hspace{-1.2pt}I\hspace{-1.2pt}I \\ \hline
    \multicolumn{1}{c|}{I\hspace{-1.2pt}I$_1$}                &
\multicolumn{1}{c|}{I\hspace{-1.2pt}I$_1$}          &
\multicolumn{1}{c|}{I\hspace{-1.2pt}I$_{\infty}$}     & \multicolumn{1}{c|}{I\hspace{-1.2pt}I$_1$}
     & \multicolumn{1}{c|}{I\hspace{-1.2pt}I$_{\infty}$} & I\hspace{-1.2pt}I\hspace{-1.2pt}I \\ \hline
    \multicolumn{1}{c|}{I\hspace{-1.2pt}I$_{\infty}$}         &
\multicolumn{1}{c|}{I\hspace{-1.2pt}I$_{\infty}$}   &
\multicolumn{1}{c|}{I\hspace{-1.2pt}I$_{\infty}$}     &
\multicolumn{1}{c|}{I\hspace{-1.2pt}I$_{\infty}$} &
\multicolumn{1}{c|}{I\hspace{-1.2pt}I$_{\infty}$} & I\hspace{-1.2pt}I\hspace{-1.2pt}I \\ \hline
    \multicolumn{1}{c|}{I\hspace{-1.2pt}I\hspace{-1.2pt}I}                 &
\multicolumn{1}{c|}{I\hspace{-1.2pt}I\hspace{-1.2pt}I}           & \multicolumn{1}{c|}{I\hspace{-1.2pt}I\hspace{-1.2pt}I}
       & \multicolumn{1}{c|}{I\hspace{-1.2pt}I\hspace{-1.2pt}I}         &
\multicolumn{1}{c|}{I\hspace{-1.2pt}I\hspace{-1.2pt}I}         & I\hspace{-1.2pt}I\hspace{-1.2pt}I
    \end{tabular}
\end{table}

\vspace{\baselineskip}

\section{Graded\ von\ Neumann\ algebra and Graded tensor product}
\vspace{\baselineskip}
In this section, we check fundamental facts about graded\ von\ Neumann\ algebra. We define some objects we study and introduce some lemmas we use later.
\vspace{\baselineskip}
\begin{dfn}\label{bo21}
Let $\mathcal{R}$ be a von Neumann algebra on a Hilbert space $\mathcal{H}$. Let 
$\Gamma$ be a self-adjoint unitary on $\mathcal{H}$ satisfies 
$\text{Ad}_{\Gamma}(\mathcal{R} )=\mathcal{R} $.\ We set 
$\theta =\text{Ad}_{\Gamma}|_\mathcal{R}$, then we say $(\mathcal{R},\
\theta)$ is a (spatially)\ graded\ von\ Neumann\ algebra with grading operator $\Gamma$.
\end{dfn}
In general, graded von Neumann algebra is a pair $(\mathcal{R},\theta)$ with $\mathcal{R}$ a von Neumann algebra and $\theta$ an involutive automorphism on $\mathcal{R}$.  
\vspace{\baselineskip}

Given a graded von Neumann algebra $(\mathcal{R},\ \theta)$,\ we set
\[\mathcal{R}^{(\sigma)}=\{x\in \mathcal{R}:\theta
(x)=(-1)^{\sigma}x\},\:\sigma=0,\:1.\] Then $\mathcal{R}$ is a direct sum of two self-adjoint $\sigma-$weakly closed linear subspaces as
$\mathcal{R}=\mathcal{R}^{(0)}\oplus \mathcal{R}^{(1)}$.\ An element of 
$\mathcal{R}^{(\sigma)}$ is said to be homogeneous of degree\ $\sigma$, or even/odd for $\sigma=0$/$\sigma=1$.
\ For a homogeneous\ $x$, its degree is denoted by $\partial x$. For  $x\in
\mathcal{R}$, we set $x^{(0)}:=\frac{x+\theta (x)}{2} $,\
$x^{(1)}:=\frac{x-\theta (x)}{2}$.\\
Let $(\mathcal{R}_1,\ \theta_1)$ and $(\mathcal{R}_2,\ \theta_2)$ be graded\
von\ Neumann\
algebras. A $*$-homomorphism\ $\varphi : \mathcal{R}_1\to \mathcal{R}_2$ is a graded homomorphism if
$\varphi (\mathcal{R}_1^{(\sigma)})\subseteq
\mathcal{R}_2^{(\sigma)}$,\ ($\sigma =0$,\ $1$).
\vspace{\baselineskip}

Balanced or Central graded von Neumann algebras are main object we study in this paper.
\vspace{\baselineskip}

\begin{dfn}\label{bo22}
Let $(\mathcal{R},\ \theta)$ be a graded\ von\ Neumann\ algebra.\
If $\mathcal{R}$ has an odd self-adjoint\ unitary,\ we say $(\mathcal{R},\ \theta)$ is balanced.\ If $Z(\mathcal{R})\cap
\mathcal{R}^{(0)}=\mathbb{C} I$ for the center $Z(\mathcal{R})$ of $\mathcal{R}$, we say
$(\mathcal{R},\ \theta)$ is central.
\end{dfn}
\vspace{\baselineskip}

Next, we define graded tensor product of graded von Neumann algebras.
\vspace{\baselineskip}

Let $(\mathcal{R}_1,\ \text{Ad}_{\Gamma_1})$\ and $(\mathcal{R}_2,\
\text{Ad}_{\Gamma_2})$ be graded von Neumann algebras acting on Hilbert spaces $\mathcal{H}_1$,\ $\mathcal{H}_2$ respectively with grading operator $\Gamma_1$, $\Gamma_2$. We define a product and involution on the algebraic tensor\ product\ $\mathcal{R}_1\odot
\mathcal{R}_2$ by
\begin{align}
(A_1\hat{\otimes }B_1)(A_2\hat{\otimes }B_2)&=(-1)^{\partial B_1
\partial  A_2}(A_1A_2\hat{\otimes }B_1B_2)\\
(A\hat{\otimes }B)^*&=(-1)^{\partial A \partial  B}A^*\hat{\otimes }B^*
\end{align}
for homogeneous simple tensors. The algebraic tensor\ product with this product and involution is denoted by $\mathcal{R}_1\hat{\odot} \mathcal{R}_2$
. This is a $*$-algebra.\ For homogeneous\ $A\in \mathcal{R}_1$,\ $B\in
\mathcal{R}_2$ we set
\[\pi(A\hat{\otimes }B)=A\Gamma_1^{\partial B}\otimes
B\] then we get faithful representation of $\mathcal{R}_1\hat{\odot} \mathcal{R}_2$ on 
$\mathcal{H}_1\otimes\mathcal{H}_2 $.\ A von Neumann algebra generated by 
$\pi(\mathcal{R}_1\hat{\odot} \mathcal{R}_2)$ is said to be a graded tensor product of $(\mathcal{R}_1,\ \mathcal{H}_1,\ \Gamma_1
)$ and $(\mathcal{R}_2,\ \mathcal{H}_2,\ \Gamma_2 )$, and denote it by $\mathcal{R}_1\hat{\otimes } \mathcal{R}_2$. $\mathcal{R}_1\hat{\otimes } \mathcal{R}_2$ is a graded von Neumann algebra with grading operator $\Gamma_1\otimes \Gamma_2$.
\vspace{\baselineskip}

The next four lemmas are Lemma A.1,\ A.2,\ A.4,\ A.5 of \cite{bobo} respectively.

\begin{lem}\label{boa1}
Let $(\mathcal{R},\ \theta)$ be a balanced\ graded\ von\ Neumann\
algebra. Assume that $\mathcal{R}$ is of type $\mu$ and $\mathcal{R}^{(0)}$ is of type\ $\lambda$, with some\ $\mu,\lambda$=I,\ I\hspace{-1.2pt}I,\ I\hspace{-1.2pt}I\hspace{-1.2pt}I, and\ both of $\mathcal{R}$, $\mathcal{R}^{(0)}$ have finite-dimensional centers. Then $\mu=\lambda$.
\end{lem}
\vspace{\baselineskip}

\begin{lem}\label{boa2}
Let $(\mathcal{R},\ \theta)$ be a central\ graded\ von\ Neumann\
algebra. Then either $Z(\mathcal{R})=\mathbb{C}I$ or $Z(\mathcal{R})$ has a self-adjoint\ unitary\
$b\in Z(\mathcal{R})\cap \mathcal{R}^{(1)}$ such that 
\[Z(\mathcal{R})\cap \mathcal{R}^{(1)}=\mathbb{C}b.\]
\end{lem}
\vspace{\baselineskip}

\begin{lem}\label{boa4}
Let $(\mathcal{R}_1,\ \text{Ad}_{\Gamma_1})$ and $(\mathcal{R}_2,\
\text{Ad}_{\Gamma_2})$ be graded von Neumann algebras acting on  $\mathcal{H}_1$,\ $\mathcal{H}_2$ respectively with grading operator $\Gamma_1$, $\Gamma_2$.\ 
Let $\mathcal{R}_1\hat{\otimes } \mathcal{R}_2$ be the graded\ tensor\
product of $(\mathcal{R}_1,\ \mathcal{H}_1,\ \Gamma_1
)$ and $(\mathcal{R}_2,\ \mathcal{H}_2,\ \Gamma_2 )$.\ Then commutant\ $(\mathcal{R}_1\hat{\otimes } \mathcal{R}_2)'$ of the graded tensor product is generated by 
\[(\mathcal{R}_1')^{(0)}\odot\mathcal{R}_2',\:(\mathcal{R}_1')^{(1)}\odot \mathcal{R}_2'\Gamma_2.\]
\end{lem}
\vspace{\baselineskip}

\begin{lem}\label{boa5}
Let $(\mathcal{R}_i,\ \text{Ad}_{\Gamma_i})$,\ $(\mathcal{L}_i,\ \text{Ad}_{W_i})$,\
$i=1$,\ $2$,\ be graded\ von\ Neumann\ algebras on $\mathcal{H}_i$,\ and $\mathcal{K}_i$ respectively.\ Let $\alpha _i:\mathcal{R}_i\to \mathcal{L}_i$,\
$i=1$,\ $2$,\ be graded\ $*$-isomorphisms.\ Suppose that $\mathcal{R}_2$ is either
(hence $ \mathcal{L}_2$ as well) balanced or trivially\ graded.\
Let $\mathcal{R}_1\hat{\otimes } \mathcal{R}_2$ be the graded\ tensor\ product of $(\mathcal{R}_1,\ \mathcal{H}_1,\ \Gamma_1 )$ and $(\mathcal{R}_2,\ \mathcal{H}_2,\ \Gamma_2 )$. Let $\mathcal{L}_1\hat{\otimes } \mathcal{L}_2$ be the graded\ tensor\ product of $(\mathcal{L}_1,\ \mathcal{K}_1,\ W_1 )$ and $(\mathcal{L}_2,\ \mathcal{K}_2, W_2 )$. Then there exists a unique $*$-isomorphism\
$\alpha _1\hat{\otimes }\alpha _2:\mathcal{R}_1\hat{\otimes }
\mathcal{R}_2\to\mathcal{L}_1\hat{\otimes } \mathcal{L}_2 $
such that 
\[(\alpha _1\hat{\otimes }\alpha _2)(A\hat{\otimes
}B)=\alpha_1(A)\hat{\otimes }\alpha _2(B)\]
for all $A\in \mathcal{R}_1$,\ $B\in \mathcal{R}_2$.
\end{lem}
\vspace{\baselineskip}

\section{Main result}
In this section, we prove main results.
\subsection{The case at least one of $\mathcal{R}_1$,\ $\mathcal{R}_2$ is of type\ I\hspace{-1.2pt}I\hspace{-1.2pt}I}

\begin{prop}\label{n3}
Let $(\mathcal{R},\ \text{Ad}_{\Gamma} )$ be a graded\ von\ Neumann\ algebra. Then
\[\mathcal{L}=\mathcal{R}^{(0)}+\mathcal{R}^{(1)}\Gamma\]is a von\
Neumann\ algebra,\ and $(\mathcal{L},\ \text{Ad}_{\Gamma})$ is a graded\ von\
Neumann\ algebra. Furthemore,
\[\mathcal{L}^{(0)}=\mathcal{R}^{(0)},\:\mathcal{L}^{(1)}=\mathcal{R}^{(1)}\Gamma.\]
\end{prop}
$<$proof$>$ It is trivial that $\mathcal{L}$ is closed under addition, scalar multiplication, and involution $*$. Since 
\[(A^{(0)}+A^{(1)}\Gamma)(B^{(0)}+B^{(1)}\Gamma)=A^{(0)}B^{(0)}-A^{(1)}B^{(1)}+(A^{(0)}B^{(1)}+A^{(1)}B^{(0)})\Gamma\in
\mathcal{L},\]
it is also closed under multiplication.\ Let $A_a^{(0)}+A_a^{(1)}\Gamma\to B$ (WOT), then we get 
$A_a^{(0)}=\frac{1}{2}(
\text{Ad}_{\Gamma}(A_a^{(0)}+A_a^{(1)}\Gamma)+A_a^{(0)}+A_a^{(1)}\Gamma)\to\frac{B+\text{Ad}_{\Gamma}(B)}{2}$ (WOT) by weak-operator continuity of $\text{Ad}_{\Gamma}$. Accordingly 
$A_a^{(1)}\Gamma\to \frac{B-\text{Ad}_{\Gamma}(B)}{2}$ (WOT) i.e. $A_a^{(1)}\to
\frac{B-\text{Ad}_{\Gamma}(B)}{2}\Gamma$.\ Since 
$\mathcal{R}^{(\sigma)}$,\ $\sigma=0$,\ $1$ is WOT closed,
\[\frac{B+\text{Ad}_{\Gamma}(B)}{2}\in \mathcal{R}^{(0)},\:
\frac{B-\text{Ad}_{\Gamma}(B)}{2}\Gamma\in \mathcal{R}^{(1)}.\] Therefore $B\in
\mathcal{L} $
 and $\mathcal{L}$ is WOT closed.\ It is clear that 
$\text{Ad}_{\Gamma}$ defines a grading of $\mathcal{L}$. Then it is also obvious that $\mathcal{L}^{(0)}=\mathcal{R}^{(0)}$,
\ $\mathcal{L}^{(1)}=\mathcal{R}^{(1)}\Gamma$. $\Box $
\vspace{\baselineskip}

\begin{prop}\label{n4}
Let $(\mathcal{R},\ \text{Ad}_\Gamma)$ be a graded\ von\ Neumann\ algebra.\
Then $\mathcal{R}$ is $*-$isomorphic to  $\mathcal{R}^{(0)}\oplus\mathcal{R}^{(1)}\Gamma$.
\end{prop}
$<$proof$>$ Since $(\mathcal{R},\ \text{Ad}_\Gamma)$ is a graded\ von\ Neumann\
algebra, $\mathcal{R} ^{(0)}\oplus \mathcal{R} ^{(1)}\Gamma$
 is a von\ Neumann\ algebra by \ref{n3}.\ Let $V=\frac{1-i}{2}I+\frac{1+i}{2}
\Gamma $, then
\begin{align}
 V^*V = VV^*
 &= (\frac{1-i}{2} I+ \frac{1+i}{2} \Gamma )( \frac{1+i}{2} I+
\frac{1-i}{2} \Gamma )\\
 &= \frac{1}{2}I-\frac{i}{2} \Gamma +\frac{i}{2} \Gamma +\frac{1}{2}I\\
 &= I
\end{align}
 so $V$ is a unitary.\ We can define a map $\mathcal{R} \twoheadrightarrow \mathcal{R}
^{(0)}\oplus \mathcal{R} ^{(1)}\Gamma $
  by $A \mapsto V^*AV$ for all $A\in \mathcal{R}$:
 \begin{align}V^*AV &= (\frac{1+i}{2}I+\frac{1-i}{2} \Gamma
)A(\frac{1-i}{2}I+\frac{1+i}{2}\Gamma )\\
  &= \frac{1}{2}A+\frac{i}{2}A\Gamma -\frac{i}{2}\Gamma
A+\frac{1}{2}\Gamma A \Gamma \\
 &= \frac{1}{2}(A+\Gamma A \Gamma )+\frac{i}{2}(A\Gamma -\Gamma A)\\
 &= A^{(0)}+iA^{(1)}\Gamma \in \mathcal{R} ^{(0)}\oplus \mathcal{R} ^{(1)}\Gamma.
 \end{align}
This map is a surjection and $*$-isomorphism.$\Box $
\vspace{\baselineskip}

\begin{prop}\label{48}
Let $(\mathcal{R} _1,\ \text{Ad}_{\Gamma _1})$ and $(\mathcal{R} _2,\ \text{Ad}_{\Gamma _2})$ be graded von Neumann algebras acting on $\mathcal{H} _1$,\ $\mathcal{H} _2$ respectively. Then 
$\mathcal{R} _1 \hat{\otimes } \mathcal{R} _2$ is $*$-isomorphic to $\mathcal{R} _2
\hat{\otimes } \mathcal{R} _1$.
\end{prop}
$<$proof$>$ Let \[V_j=\frac{1-i}{2}I+\frac{1+i}{2}\Gamma_j ,\  (j=1,\:2)\]
\[V=\frac{1-i}{2}I+\frac{1+i}{2}\Gamma_1 \otimes \Gamma_2,\] $U$ be a unitary $\mathcal{H} _1\otimes\mathcal{H} _2\to \mathcal{H} _2\otimes\mathcal{H} _1$ such that
\[U(x \otimes y)=y\otimes x,\ (x\in\mathcal{H} _1,\:y\in\mathcal{H}
_2).\]
The equation \[\phi=\text{Ad}_U\text{Ad}_{\Gamma_1 \otimes \Gamma_2}\text{Ad}_{V}\text{Ad}_{V_1\otimes
V_2}\] defines a map $\phi:\mathcal{R} _1 \hat{\otimes } \mathcal{R} _2  \to \mathcal{R} _2
\hat{\otimes } \mathcal{R} _1$. We determine the image of $\mathcal{R} _1 \hat{\otimes } \mathcal{R} _2$ under $\phi$.\ Let $A_\sigma\in \mathcal{R} _1^{(\sigma)}$,\ $A_\sigma\in
\mathcal{R} _2^{(\sigma)}$,\ $\sigma=0$,\ $1$. A Straightforward calculation shows that 
\[\phi(A_0\otimes B_0)=B_0\otimes A_0,\:\phi(A_0 \Gamma_1\otimes
B_1)=B_1\otimes A_0\]
\[\phi(A_1\otimes B_0)=B_0 \Gamma_2\otimes
A_1,\:\phi(A_1\Gamma_1\otimes B_1)=-B_1 \Gamma_2\otimes A_1.\]
Since $\phi$ is a $*$-isomorphism, $\phi(\mathcal{R} _1 \hat{\otimes } \mathcal{R} _2)$ is a von Neumann algebra which is generated by
\[\mathcal{R} _2\odot \mathcal{R} _1^{(0)},\:\mathcal{R} _2
\Gamma_2\odot \mathcal{R} _1^{(1)}\] i.e. $\mathcal{R} _2 \hat{\otimes } \mathcal{R} _1$.$\Box $
\vspace{\baselineskip}

\begin{prop}\label{45}
Let $(\mathcal{R} _1,\ \text{Ad}_{\Gamma _1})$ and $(\mathcal{R} _2,\ \text{Ad}_{\Gamma _2})$ be graded\ von\ Neumann\ algebra acting on $\mathcal{H} _1$,\ $\mathcal{H} _2$ respectively.\ Assume that 
$\mathcal{R}$ is $(\mathcal{R}_1\hat{\otimes}\mathcal{R}_2)^{(0)}\oplus(\mathcal{R}_1\hat{\otimes}\mathcal{R}_2)^{(1)}(\Gamma_1\otimes\Gamma_2)$
 and $\mathcal{S} $ is a von\ Neumann\ algebra generated by $I\odot \mathcal{R} _2^{(0)}$ and $I\odot \mathcal{R}
_2^{(1)}\Gamma _2$. Then there is a family $\{\Phi  _z:z \in \mathcal{H} _1,\ \|z\|=1\}$
 of conditional\ expectations\ $\Phi  _z:\mathcal{R} \twoheadrightarrow
\mathcal{S} $ which satisfies following $2$ properties.\\
(i)\ For all $z$, $\Phi _z$ is weak-operator continuous on the unit ball $(\mathcal{R} )_1$,\\
(ii)\ If $T \in \mathcal{R} ^+ $and $T \neq 0$, then $\Phi _z(T) \neq 0$ for some $z$.
\end{prop}
$<$proof$>$ The proof is similar to that of 11.2.24\ PROPOSITION in \cite{krkr}.\ Assume $z \in
\mathcal{H} _1$,\ $\|z\|=1$, and 
\[E(T)=\frac{T+\text{Ad}_{\Gamma _1\otimes I}T}{2}\]
for $T \in \mathcal{R} _1$. The equation 
\[b_{zT}(x,\ y)=\langle E(T)(z\otimes x),\ z\otimes y\rangle\ (x,\ y \in
\mathcal{H} _2)\] defines a conjugate\ bilinear\ functional\ $b_{zT}$ on $\mathcal{H}_2$. Since $b_{zT}$ is bounded, $b_{zT}$ corresponds to a bounded\
linear\ map\ $\Psi  _z(T)$ on $\mathcal{H} _2$ which satisfies  $b_{zT}(x,\ y)=\langle\Psi  _z(T)x,\ y\rangle$ for all $x$, $y\in \mathcal{H} _2$.\ It is obvious that 
$\Psi _z:\mathcal{R} \to \mathfrak{B} (\mathcal{H} _2)$ is weak-operator continuous, positive,\ and
$\Psi _z(I)=I$.\
It follows that $\mathcal{S}\subseteq\mathcal{R}$, since 
\[I\otimes (A_2+ B_2\Gamma_2)=I\otimes A_2+(\Gamma_1\otimes
B_2)(\Gamma_1\otimes\Gamma_2)\] for all $A_2\in\mathcal{R}_2^{(0)}$,\
$B_2\in \mathcal{R} _2^{(1)}$. Let 
\[\mathcal{M}=\mathcal{R}_2^{(0)}\oplus\mathcal{R}_2^{(1)}\Gamma_2.\]\
We show that 
\begin{align}
\Psi_z(T)\in\mathcal{M},\ \Psi_z((I\otimes X)T(I\otimes Y))=X\Psi_z(T)Y \label{tonma}
\end{align}
for all $T \in \mathcal{R}  $,\ $X$,\ $Y \in \mathcal{M} $.\ Since $\Psi_z$ is weak-operator continuous, it will suffice to show that in the case $T$ has the form 
$A\otimes(B+B'\Gamma_2)$,\ $C\Gamma_1\otimes D\Gamma_2$,\
$E\Gamma_1\otimes F$,\ $A\in \mathcal{R}_1^{(0)}$,\ $C,\
E\in\mathcal{R}_1^{(1)}$
,\ $B,\ D \in \mathcal{R} _2^{(0)}$,\ $B',\ F \in \mathcal{R}
_2^{(1)}$. First, let $R \in \mathfrak{B} (\mathcal{H} _1)$,\
$S \in \mathfrak{B} (\mathcal{H} _2)$, then $Z=\langle Rz,z\rangle S$ satisfies 
\[\langle Zx,\ y\rangle =\langle (R\otimes S)(z\otimes x),\ (z\otimes y)\rangle=\langle Rz,\ z\rangle\langle Sx,\ y\rangle\] for all $x\in \mathcal{H}_1$, $y\in \mathcal{H}_2$. Conversely, from the uniqueness of the Riesz representation theorem, such $Z$ is equal to $\langle Rz,z\rangle S$. Therefore  
\[\Psi_z(A\otimes(B+B'\Gamma_2)+C\Gamma_1\otimes D\Gamma_2+E\Gamma
_1\otimes F)=\langle Az,\ z\rangle (B+B'\Gamma_2).\]
Suppose $X,\ Y\in\mathcal{M}$, then 
\begin{align}
&\Psi _z((I\otimes X)(A\otimes (B+B'\Gamma _2)+C\Gamma_1\otimes
D\Gamma _2+E\Gamma _1\otimes F)(I\otimes Y))\\
&=\Psi _z((I\otimes X)(A\otimes (B+B'\Gamma _2))(I\otimes Y))\\
&=\langle Az,\ z\rangle X (B+B'\Gamma _2)Y\\
&=X\Psi _z(A\otimes (B+B'\Gamma _2))Y
\end{align}
It follows \eqref{tonma}. The equation 
\[\Phi _z(T)=I\otimes \Psi _z(T)\]
defines a map $\mathcal{R}\to \mathcal{S} $ such that $\Phi _z(I)=I$ and 
$\Phi _z(T)\in \mathcal{S} $ because $\Psi _z(I)=I$ and $\Psi_z(T)\in\mathcal{M}$.\ Since the map $A\mapsto I\otimes A$
,\ $\mathcal{M} \to \mathcal{R} $ is weak-operator continuous on the unit ball $(\mathcal{M})_1$, $\Phi _z$ is a conditional\
expectation which is weak-operator continuous on $(\mathcal{R} )_1$.\\
Finally, we assume that $T \in (\mathcal{R} )^+$ and $T \neq 0$. Since 
\[E(T)=\frac{T+\text{Ad}_{\Gamma _1\otimes I}T}{2}\] is positive and $\neq 0$, there are vectors  $z \in \mathcal{H} _1$,\ $x \in \mathcal{H}_2$ such that $E(T)(z\otimes x)\neq 0$. Then
\begin{align}
0\neq & \|E(T)^{1/2}(z\otimes x)\|\\
& = \langle E(T)(z\otimes x),\ z\otimes x\rangle\\
& = \langle \Psi _z(T)x,\ x\rangle.
\end{align}
It follows that $\Phi _z(T)=I\otimes \Psi _z(T) \neq 0$.$\Box $
\vspace{\baselineskip}

\begin{prop}\label{46}
If $(\mathcal{R} _1,\ \text{Ad}_{\Gamma _1})$ and $(\mathcal{R} _2,\ \text{Ad}_{\Gamma
_2})$ are graded\ von\ Neumann\ algebras and $\mathcal{R}_2$ is of type I\hspace{-1.2pt}I\hspace{-1.2pt}I, then $\mathcal{R} _1 \hat{\otimes }
\mathcal{R} _2$ is of type\ I\hspace{-1.2pt}I\hspace{-1.2pt}I.
\end{prop}
$<$proof$>$ We adopt the notation of \ref{45}. By using 11.2.25\ PROPOSITION in \cite{krkr}, if $\mathcal{S}$ is of 
type\ I\hspace{-1.2pt}I\hspace{-1.2pt}I, then it follows that  $\mathcal{R} $  in \ref{45} is of type\ I\hspace{-1.2pt}I\hspace{-1.2pt}I . Note that the existence of the family of conditional\ expectations which satisfies (i),\ (ii) is guaranteed by \ref{45}, and that the fact that $\mathcal{S} $ is of type\
I\hspace{-1.2pt}I\hspace{-1.2pt}I is guaranteed by $\mathcal{R} _2$ is of type\ I\hspace{-1.2pt}I\hspace{-1.2pt}I because 
\[\mathcal{S} \cong I \overline{\otimes } \mathcal{M} \cong
\mathcal{R} _2 \]by \ref{n4} .\ It follows that $\mathcal{R} _1 \hat{\otimes } \mathcal{R} _2$ is $*$-isomorphic to $\mathcal{R}
$ by \ref{n4}.$\Box $
\vspace{\baselineskip}

\begin{prop}\label{32}
Let $(\mathcal{R} _1,\ \text{Ad}_{\Gamma _1})$ and $(\mathcal{R} _2,\ \text{Ad}_{\Gamma
_2})$ be graded\ von\ Neumann\ algebras.\ Then $\text{Ad}_{I\otimes \Gamma _2}$ defines a grading on $\mathcal{R} _1 \hat{\otimes }
\mathcal{R} _2$. If $(\mathcal{R} _1 \hat{\otimes } \mathcal{R} _2)^{<0>}$ and $(\mathcal{R}_1 \hat{\otimes } \mathcal{R} _2)^{<1>}$ denote the even/odd elements defined by this grading\ respectively, then
\[(\mathcal{R} _1 \hat{\otimes } \mathcal{R} _2)^{<0>}=\mathcal{R}
_1\overline{\otimes }(\mathcal{R} _2)^{(0)} \]
holds. If $(\mathcal{R} _2,\ \text{Ad}_{\Gamma
_2})$ is balanced, then $(\mathcal{R} _1 \hat{\otimes }
\mathcal{R} _2, \text{Ad}_{I\otimes \Gamma _2})$ is also balanced.
\end{prop}
$<$proof$>$ We consider the graded von Neumann algebra $(\mathcal{R} _1 \hat{\otimes }
\mathcal{R} _2, \text{Ad}_{I\otimes \Gamma _2})$.\ It is obvious that 
$\mathcal{R} _1\overline{\otimes }(\mathcal{R} _2)^{(0)} \subseteq
(\mathcal{R} _1 \hat{\otimes } \mathcal{R} _2)^{<0>}$.\ Each element in 
$(\mathcal{R} _1 \hat{\otimes } \mathcal{R} _2)^{<0>}$ can be approximated by sum of simple\ tensors in weak-operator topology.\ Let 
$\sum A_i\Gamma _1^{\partial B_i}\otimes B_i$ be sum of simple\ tensors, then 
\[\frac{\text{id}+\text{Ad}_{I\otimes \Gamma _2}}{2} (\sum A_i\Gamma _1^{\partial
B_i}\otimes B_i)=\sum _{\partial B_i=0}A_i\otimes B_i \in
\mathcal{R} _1\overline{\otimes }(\mathcal{R} _2)^{(0)}
\]where sum of right-hand side is taken only for those satisfying $\partial B_i=0$.\\
Let $(\mathcal{R} _2,\ \text{Ad}_{\Gamma
_2})$ be balanced graded von Neumann algebra, and let
$U_2\in \mathcal{R} _2^{(1)}$ be a self-adjoint\ unitary in $\mathcal{R}_2$. Then
$\Gamma_1\otimes U_2$ is odd self-adjoint\ unitary in $(\mathcal{R}_1\hat{\otimes}\mathcal{R}_2,\
\text{Ad}_{I\otimes \Gamma _2})$. Consequently $(\mathcal{R}_1\hat{\otimes}\mathcal{R}_2,\
\text{Ad}_{I\otimes \Gamma _2})$ is balanced.
$\Box$
\vspace{\baselineskip}

\begin{prop}\label{31}
Let $(\mathcal{R} _1,\ \text{Ad}_{\Gamma _1})$ and $(\mathcal{R} _2,\ \text{Ad}_{\Gamma
_2})$ be graded\ von\ Neumann\ algebras.
If $\mathcal{R} _1\overline{\otimes }(\mathcal{R} _2)^{(0)}$ is of type\ I\hspace{-1.2pt}I\hspace{-1.2pt}I, then $\mathcal{R} _1 \hat{\otimes } \mathcal{R} _2$ is also of type\ I\hspace{-1.2pt}I\hspace{-1.2pt}I.
\end{prop}
$<$proof$>$ Suppose that $E\in\mathcal{R}_1 \hat{\otimes}\mathcal{R}
_2$ is a non-zero projection which is finite relative to $\mathcal{R}
_1 \hat{\otimes } \mathcal{R} _2$. Since $\text{Ad}_{I\otimes \Gamma _2}$ is a $*$-isomorphism, $\text{Ad}_{I\otimes \Gamma _2}(E)$ is also finite relative to $\mathcal{R} _1 \hat{\otimes }\mathcal{R} _2$. Now,
\[\text{Ad}_{I\otimes \Gamma _2}(E\lor \text{Ad}_{I\otimes \Gamma
_2}(E))=\text{Ad}_{I\otimes \Gamma _2}(E)\lor
\text{Ad}_{I\otimes \Gamma _2}(\text{Ad}_{I\otimes \Gamma _2}(E))=E\lor \text{Ad}_{I\otimes
\Gamma _2}(E)\]
so that $E\lor \text{Ad}_{I\otimes \Gamma _2}(E)$ is a element of $\mathcal{R}
_1\overline{\otimes }(\mathcal{R} _2)^{(0)}$ from \ref{32}.\ From
\cite{krkr} 6.3.8. THEOREM, $E\lor \text{Ad}_{I\otimes \Gamma
_2}(E)$ is finite relative to $\mathcal{R} _1 \hat{\otimes } \mathcal{R} _2$, therefore it is also finite relative to 
$\mathcal{R} _1\overline{\otimes }(\mathcal{R} _2)^{(0)}$.\ This contradicts the fact that $\mathcal{R} _1\overline{\otimes
}(\mathcal{R} _2)^{(0)}$ is of type\ I\hspace{-1.2pt}I\hspace{-1.2pt}I.$\Box $
\vspace{\baselineskip}

\begin{prop}\label{47}
Let $(\mathcal{R} _1,\ \text{Ad}_{\Gamma _1})$ and $(\mathcal{R} _2,\ \text{Ad}_{\Gamma
_2})$ be graded\ von\ Neumann\ algebras. If 
$\mathcal{R}_1$ is of type\ I\hspace{-1.2pt}I\hspace{-1.2pt}I, then $\mathcal{R} _1 \hat{\otimes }
\mathcal{R} _2$ is also of type\ I\hspace{-1.2pt}I\hspace{-1.2pt}I.
\end{prop}
$<$proof$>$ If $\mathcal{R} _1$ is of type\
I\hspace{-1.2pt}I\hspace{-1.2pt}I, $\mathcal{R}_1\overline{\otimes} \mathcal{R}_2^{(0)}$ is also of type\ I\hspace{-1.2pt}I\hspace{-1.2pt}I from \cite{krkr} 11.2.26. PROPOSITION. Therefore $\mathcal{R} _1 \hat{\otimes }
\mathcal{R} _2$ is also of type\ I\hspace{-1.2pt}I\hspace{-1.2pt}I from  \ref{31}.$\Box $
\vspace{\baselineskip}

From \ref{46} and \ref{47}, in the case at least one of $\mathcal{R}_1$,\ $\mathcal{R}_2$ is of type\ I\hspace{-1.2pt}I\hspace{-1.2pt}I, 
$\mathcal{R}_1 \hat{\otimes} \mathcal{R}_2$ is of type\
I\hspace{-1.2pt}I\hspace{-1.2pt}I. We can also show this by using \ref{48}.

\subsection{The case both $(\mathcal{R} _1,\ \text{Ad}_{\Gamma_1})$ and $(\mathcal{R}_2,\
\text{Ad}_{\Gamma_2})$ are central and both $\mathcal{R}_1$,\ $\mathcal{R}_2$ are not factors}

\begin{prop}\label{13}
Let $(\mathcal{R} ,\ \text{Ad}_{\Gamma }) $ be a graded\ von\ Neumann\
algebra. Then $(Z(\mathcal{R}) ,\ \text{Ad}_{\Gamma }) $ is also a graded\ von\
Neumann\ algebra.
\end{prop}
$<$proof$>$ We show that $Z(\mathcal{R})$ is closed under $\text{Ad}_{\Gamma }$.\ Assume that  $A \in
Z(\mathcal{R})$. Then for each $B \in \mathcal{R}$
\[B\Gamma A \Gamma = \Gamma (B^{(0)}-B^{(1)}) A \Gamma = \Gamma A
(B^{(0)}-B^{(1)})  \Gamma = \Gamma A \Gamma B.\]
$\Box $
\vspace{\baselineskip}

\begin{prop}\label{16}
Let $(\mathcal{R} ,\  \text{Ad}_\Gamma )$ be a graded\ von\ Neumann\
algebra such that $(Z(\mathcal{R}) ,\  \text{Ad}_\Gamma )$ is trivially\ graded. Then $Z(\mathcal{R} )\subseteq Z(\mathcal{R}^{(0)} )$.
\end{prop}
$<$proof$>$
\[Z(\mathcal{R})\subseteq \mathcal{R}^{(0)}\cap
\mathcal{R}^{'(0)}\subseteq\mathcal{R}^{(0)}\cap
\mathcal{R}^{(0)'}=Z(\mathcal{R}^{(0)}).\]
$\Box $
\vspace{\baselineskip}

\begin{prop}\label{17}
Let $(\mathcal{R} ,\  \text{Ad}_\Gamma )$ be a graded\ von\ Neumann\
algebra such that $(Z(\mathcal{R}) ,\  \text{Ad}_\Gamma )$ is balanced. Then
$Z(\mathcal{R}^{(0)} )\subseteq Z(\mathcal{R} )$.
\end{prop}
$<$proof$>$ Assume that $A \in Z(\mathcal{R}^{(0)} ) $.\ Since $(Z(\mathcal{R}) ,\
\text{Ad}_\Gamma )$ is balanced, there is a self-adjoint\ unitary $ U \in Z(\mathcal{R} )^{(1)}$. For each $B \in \mathcal{R} ^{(1)}$, $B = BUU$,\
$BU \in \mathcal{R}^{(0)}$, so that
\[AB = BUAU= BA.\]Therefore $A \in \mathcal{R} ^{(1)'}$, and $A \in
\mathcal{R} '$. Then $A \in Z(\mathcal{R} )$.
$\Box $
\vspace{\baselineskip}

\begin{prop}\label{19}
Let $(\mathcal{R} ,\  \text{Ad}_\Gamma )$ be a graded\ von\ Neumann\ algebra such that $(Z(\mathcal{R}) ,\  \text{Ad}_\Gamma )$ is balanced. Then abelian\ projections in 
$\mathcal{R} ^{(0)}$ are also abelian in $\mathcal{R} $.
\end{prop}
$<$proof$>$ Let $E\in \mathcal{R} ^{(0)}$ be an abelian\
projection in $\mathcal{R} ^{(0)}$,\ and let $U\in Z(\mathcal{R})^{(1)}$ be a self-adjoint\ unitary.\ Then for each $A$,\ $B\in \mathcal{R}$
\[EAEBE=(EA^{(0)}E+EA^{(1)}UEU)(EB^{(0)}E+EB^{(1)}UEU).\]
Since $E$ is an abelian\ projection in $\mathcal{R} ^{(0)}$ and $A^{(1)}U$,\
$B^{(1)}U\in\mathcal{R} ^{(0)}$,\ $U\in
Z(\mathcal{R})$, it follows that $EA^{(0)}E$ commute with $EB^{(0)}E+EB^{(1)}UEU$ and 
$EA^{(1)}UEU$ commute with $EB^{(0)}E+EB^{(1)}UEU$
. Thus $EAEBE=EBEAE$, and $E$ is an ableian projection in $\mathcal{R} $.$\Box $
\vspace{\baselineskip}

\begin{prop}\label{11}
Let $(\mathcal{R} ,\ \text{Ad}_{\Gamma }) $ be a central\ graded\ von\ Neumann\
algebra. Assume that $\mathcal{R}$ is of type\ I$_n$($n$ is a cardinal number) and not a factor. Then $\mathcal{R}^{(0)}$ is a type\ I$_n$\ factor. Furthermore, ablelian\ projections in $\mathcal{R}^{(0)}$ are also abelian in $\mathcal{R}$.
\end{prop}
$<$proof$>$ The proof of the first argument is similar to the proof of \cite{bobo}Prop.2.9. First, we shall see that $\mathcal{R}^{(0)}$ is a factor.\ 
Since $(\mathcal{R} ,\ \text{Ad}_{\Gamma }) $ is central and $(Z(\mathcal{R}) ,\ \text{Ad}_\Gamma )$ is balanced, 
\[Z(\mathcal{R} ^{(0)})\subseteq \mathcal{R} ^{(0)}\cap
Z(\mathcal{R})=\mathbb{C} I\]from \ref{17}.
Since the reverse inclusion is apparent, it follows that $Z(\mathcal{R}^{(0)})=\mathbb{C}
I$ and $\mathcal{R}^{(0)}$ is a factor.\ Since $\mathcal{R}$ is of type\
I and balanced,\ $\mathcal{R}^{(0)}$ is a
factor (and especially of type\ I or I\hspace{-1.2pt}I or\ I\hspace{-1.2pt}I\hspace{-1.2pt}I),\ and both
$\mathcal{R}$, $\mathcal{R}^{(0)}$ have finite-dimensional centers, it follows that 
$\mathcal{R}^{(0)}$ is a type\ I factor from \cite{bobo} Lem.A.1.\ Assume that $\mathcal{R}^{(0)}$ is of type\
I$_m$($m$ is a cardinal number).\ Since $\mathcal{R}^{(0)}$ is of type\ I$_m$, there is an orthogonal family $\{E_a:a\in \mathbb{A} \}$ of mutually equivalent ablelian\
projection in $\mathcal{R}^{(0)}$ such that $|\mathbb{A} |=m$ and $\sum_{a} E_a =I$.\
$\{E_a\}$ is mutually orthogonal in $\mathcal{R}^{(0)}$, so also in $\mathcal{R}$. 
Since $(Z(\mathcal{R}) ,\  \text{Ad}_\Gamma
)$ is balanced, each $E_a$ is also abelian in $\mathcal{R}$ from \ref{19}. Therefore $I$ is written in the form $\sum_{a} E_{a} =I$ by mutually orthogonal family $\{E_a:a\in
\mathbb{A} \}$,\ $|\mathbb{A} |=m$ of abelian projections in 
$\mathcal{R}$. Thus $\mathcal{R}$ is also of type type\ I$_m$ and $m=n$
from \cite{krkr} 6.5.2. THEOREM.
We have already shown that ablelian\ projections in $\mathcal{R}^{(0)}$ are also abelian in $\mathcal{R}$.$\Box$
\vspace{\baselineskip}

\begin{prop}\label{12}
Let $(\mathcal{R} _1,\ \text{Ad}_{\Gamma _1}) $ and $(\mathcal{R} _2,\ \text{Ad}_{\Gamma
_2}) $ be central\ graded\ von\ Neumann\ algebras acting on $\mathcal{H} _1$,\ $\mathcal{H} _2$ respectively. If $\mathcal{R} _1$ is of type\
I$_m$($m$ is a cardinal number) and not a factor, and $\mathcal{R} _2$ is of type\ I$_n$($n$ is a cadinal number) and not a factor, then $\mathcal{R} _1\hat{\otimes } \mathcal{R} _2$ is of type\ I$_{2mn}$.
\end{prop}
$<$proof$>$ Since both $(\mathcal{R} _1,\ \text{Ad}_{\Gamma _1}) $ and $(\mathcal{R} _2,\
\text{Ad}_{\Gamma _2}) $ are central\ graded\ von\ Neumann\ algebras, and 
$\mathcal{R}_1$ and $\mathcal{R}_2$ are of type\ I$_m$,\ type\
I$_n$ respectively and not a factor, $\mathcal{R} _1^{(0)}$ and  $\mathcal{R} _2^{(0)}$
 are of type\ I$_m$,\ type\ I$_n$\ factor respectively from \ref{11}.\ There is a family $\{E_a\}_{a\in \mathbb{A}}\subseteq
\mathcal{R} _1^{(0)}$ of abelian projections equivalent in $\mathcal{R} _1^{(0)}$ such that $\sum_{a} E_{a}
=I$, $|\mathbb{A}|=m$.\ And there is a similar family $\{F_b\}_{b\in\mathbb{B}}\subseteq \mathcal{R} _2^{(0)}$ such that $|\mathbb{B}|=n$ for $\mathcal{R} _2^{(0)}$. Suppose that
\ $U_1 \in Z(\mathcal{R} _1)\cap \mathcal{R} _1^{(1)}$ is a self-adjoint\ unitary. Let
\[G_{a,\ b} = \frac{1+U_1}{2}E_a\otimes F_b,\:G_{a,\ b}' =
\frac{1-U_1}{2}E_a\otimes F_b, \]
then $G_{a,\ b}$ and $G_{a,\ b}'$ are abelian\ projections in $\mathcal{R} _1\hat{\otimes }
\mathcal{R} _2$.\ We shall show $G_{a,\ b}$ is an abelian\ projection.\ For all $A \in \mathcal{R} _1$,\ $B\in \mathcal{R} _2^{(1)}$,
\[G_{a,\ b}(A\Gamma_1 \otimes B)G_{a,\ b}=\frac{1+U_1}{2}E_a A
\Gamma_1 \frac{1+U_1}{2}E_a \otimes F_b B F_b=0\] so that
\[G_{a,\ b}(\mathcal{R} _1\hat{\otimes } \mathcal{R} _2)G_{a,\ b}
=G_{a,\ b}(\mathcal{R} _1\hat{\otimes } \mathcal{R} _2^{(0)})G_{a,\ b}.\]
In fact,\ $G_{a,\ b}(\mathcal{R} _1\hat{\otimes } \mathcal{R} _2)G_{a,\ b}
\supseteq G_{a,\ b}(\mathcal{R} _1\hat{\otimes } \mathcal{R} _2^{(0)})G_{a,\ b}$ is obvious. We show the reverse inclusion as follows. Suppose $X \in \mathcal{R} _1\hat{\otimes } \mathcal{R} _2$,\ and let $\{X_t\}$ be a net of the elements which is the sum of elements of the form $A \Gamma_1^{\partial B}\otimes B$ which weak-operator converges to $X$, then
\[G_{a,\ b}X G_{a,\ b} = \lim _{t}G_{a,\ b}X_{t}G_{a,\ b}.\]
Since $G_{a,\ b}X_{t}G_{a,\ b} \in G_{a,\ b}(\mathcal{R} _1\hat{\otimes }
\mathcal{R} _2^{(0)})G_{a,\ b}$ and 
$G_{a,\ b}(\mathcal{R} _1\hat{\otimes } \mathcal{R} _2^{(0)})G_{a,\
b}$ is weak-operator closed in $\mathfrak{B} (\mathcal{H} )$, this completes the proof of the reverse inclusion.\\
 Since $(\mathcal{R} _1,\ \text{Ad}_{\Gamma _1}) $ is a central\ graded\ von\ Neumann\ algebra and
$\mathcal{R}_1$ is of type\ I$_m$ and not a factor,  all of $E_a$ are abelian\ projections of 
$\mathcal{R} _1$ by \ref{11}. Thus 
$(E_a \otimes F_b)\mathcal{R} _1\hat{\otimes } \mathcal{R}
_2^{(0)}(E_a \otimes F_b)$ is abelian. We note that
$G_{a,\ b} \in (Z(\mathcal{R} _1)\hat{\otimes } Z(\mathcal{R}
_2^{(0)}))(E_a\otimes F_b)$. Hence 
\[G_{a,\ b}(E_a \otimes F_b)(\mathcal{R} _1\hat{\otimes } \mathcal{R}
_2^{(0)})(E_a \otimes F_b)G_{a,\ b}
=G_{a,\ b}(\mathcal{R} _1\hat{\otimes } \mathcal{R} _2^{(0)})G_{a,\ b}\]
is abelian.\ $G_{a,\ b}'$ is also abelian in the same way.\ Since $\displaystyle\sum_{a,\ b}G_{a,\ b}
+\sum_{a,\ b}G_{a,\ b}'=I$, all we have to do is showing they are mutually equivalent. To see this, it suffices to show that the central supports of each $G_{a,\ b}$,\ $G_{a,\
b}'$ are $I$ in $\mathcal{R} _1\hat{\otimes } \mathcal{R} _2$ from \cite{krkr} 6.4.6PROPOSITON.\\
It suffices to show that $\overline{(\mathcal{R} _1\hat{\otimes }
\mathcal{R} _2)G_{a,\ b}
(\mathcal{H} _1 \otimes \mathcal{H} _2)} = \mathcal{H} _1 \otimes
\mathcal{H} _2$ from \cite{krkr} 5.5.2.PROPOSITION to see that the central support of $G_{a,\ b}$ is $I$. We shall show that $x\otimes
y\in C_{G_{a,\ b}}(\mathcal{H}_1\otimes\mathcal{H}_2)$ for each $x\in\mathcal{H}_1$,\ $y\in\mathcal{H}_2$.\ Since $\mathcal{R} _1^{(0)}$ is a factor, the central\ support of $E_a$ in $\mathcal{R}
_1^{(0)}$ is $I$. Thus $\overline{\mathcal{R} _1^{(0)}E_a (\mathcal{H} _1)}
=\mathcal{H} _1$ from \cite{krkr} 
5.5.2.PROPOSITION. Given a positive $\epsilon >0$, there are 
$x_1 \in \mathcal{H} _1$,\ $A \in \mathcal{R} _1^{(0)}$ such that $\|AE_a
x_1 -x\|<\epsilon $. Similarly, there are $y_1,\ y_2 \in \mathcal{H} _2$
,\ $B_1,\ B_2 \in \mathcal{R} _2^{(0)}$ such that $\|B_1 F_b y_1 -y \|<\epsilon $,\ $\|B_2 F_b y_2 -U_2 y\|<\epsilon $, where $U_2 \in Z(\mathcal{R}
_2)\cap \mathcal{R} _2^{(1)}$ is a self-adjoint\ unitary. Note that
\begin{align}
&\left\|A\dfrac{1+U_1}{2}E_a x_1 \otimes B_1 F_b y_1 - A\frac{1+U_1}{2}E_a
x_1 \otimes y\right\|\\
&\leq \|A\frac{1+U_1}{2}E_a x_1 \|\|B_1 F_b y_1 - y\|\\
&\leq \|A\frac{1+U_1}{2}E_a x_1\| \epsilon=\|\frac{1+U_1}{2}AE_a x_1\|
\epsilon \\
&\leq \|AE_a x_1\| \epsilon \leq (\|x\|+\epsilon )\epsilon. \label{aho}
\end{align}
On the other hand, 
\[(A\Gamma _1 \otimes U_2 B_2)G_{a,\ b}(\Gamma _1 x_1 \otimes y_2)
 = A\frac{1-U_1}{2}E_a x_1 \otimes U_2 B_2 F_b y_2\]thus
\begin{align}
&\|(A\Gamma _1 \otimes U_2 B_2)G_{a,\ b}(\Gamma x_1 \otimes y_2)
- A\frac{1-U_1}{2}E_a x_1 \otimes y \|\\
&\leq \| A\frac{1-U_1}{2}E_a x_1\|\|U_2 B_2 F_b y_2 - y\|\\
&\leq \|A\frac{1-U_1}{2}E_a x_1\|\epsilon \leq (\|x\|+\epsilon
)\epsilon.\label{baka}
\end{align}
Since $\epsilon >0$ is arbitrary,
\[ A\frac{1+ U_1}{2}E_a x_1 \otimes y \in C_{G_{a,\ b}}(\mathcal{H}
_1\otimes \mathcal{H} _2) \]
\[ A\frac{1- U_1}{2}E_a x_1 \otimes y \in C_{G_{a,\ b}}(\mathcal{H}
_1\otimes \mathcal{H} _2) \]from \eqref{aho},\ \eqref{baka} respectively . Consequently
\[x\otimes y \in C_{G_{a,\ b}}(\mathcal{H} _1\otimes \mathcal{H} _2) \]and
\[ C_{G_{a,\ b}}(\mathcal{H} _1\otimes \mathcal{H} _2) = \mathcal{H}
_1\otimes \mathcal{H} _2 .\]
\ Similarly the central support of $G_{a,\ b}'$ is also $I$. Thus $\mathcal{R} _1\hat{\otimes } \mathcal{R} _2$ is of type\
I$_{2mn}$.$\Box$
\vspace{\baselineskip}

\begin{prop}\label{18}
Let $(\mathcal{R} _1,\ \text{Ad}_{\Gamma _1})$ and $(\mathcal{R} _2,\ \text{Ad}_{\Gamma
_2})$ be graded\ von Neumann\ algebras. If both
$(Z(\mathcal{R} _1),\ \text{Ad}_{\Gamma _1})$ and $(Z(\mathcal{R} _2),\
\text{Ad}_{\Gamma _2})$ are balanced, then 
$Z(\mathcal{R} _1\hat{\otimes } \mathcal{R} _2) = Z(\mathcal{R}
_1)^{(0)}\overline{\otimes }Z(\mathcal{R} _2)^{(0)} $.
\end{prop}
$<$proof$>$ It is always true that $Z(\mathcal{R} _1)^{(0)}\overline{\otimes }Z(\mathcal{R}
_2)^{(0)} \subseteq Z(\mathcal{R} _1\hat{\otimes } \mathcal{R} _2)$ because each simple\ tensors in left-hand side are element of the right-hand side. To prove the reverse inclusion, it suffices to show that 
$(Z(\mathcal{R} _1)^{(0)}\overline{\otimes }Z(\mathcal{R}
_2)^{(0)})' \subseteq Z(\mathcal{R} _1\hat{\otimes } \mathcal{R} _2)'$
from the double commutant theorem (\cite{krkr} 5.3.1. THEOREM). The commutant of the left-hand side is 
$Z(\mathcal{R} _1)^{(0)'}  \overline{\otimes } Z(\mathcal{R} _2)^{(0)'} $ from \cite{krkr} 11.2.16THEOREM. Note that
$Z(\mathcal{R} _1)^{(0)} = \mathcal{R} _1 \cap \mathcal{R} _1 '\cap
\{\Gamma _1 \}'$ hence $Z(\mathcal{R} _1)^{(0)'}$ is generated by
$\mathcal{R} _1 \cup  \mathcal{R} _1 '\cup  \{\Gamma _1 \}''$.\ Similarly 
$Z(\mathcal{R} _2)^{(0)'}$ is generated by $\mathcal{R} _2 \cup  \mathcal{R} _2 '\cup  \{\Gamma _2 \}''$. On the other hand, 
$Z(\mathcal{R} _1\hat{\otimes } \mathcal{R} _2)'$ is generated by
$\mathcal{R} _1\hat{\otimes } \mathcal{R} _2\cup (\mathcal{R}
_1\hat{\otimes } \mathcal{R} _2)'$, 
$\mathcal{R} _1\hat{\otimes } \mathcal{R} _2$ is generated by 
$\mathcal{R} _1\odot \mathcal{R} _2^{(0)}\cup \mathcal{R} _1\Gamma
_1\odot \mathcal{R} _2^{(1)}$, and 
$(\mathcal{R} _1\hat{\otimes } \mathcal{R} _2)'$ is generated by 
$(\mathcal{R} _1 ')^{(0)}\odot \mathcal{R} _2 '\cup (\mathcal{R} _1
')^{(1)}\odot R_2 ' \Gamma _2$ from \cite{bobo} LemmaA.4.
Thus $Z(\mathcal{R} _1\hat{\otimes } \mathcal{R} _2)'$ is generated by 
$\mathcal{R} _1\odot \mathcal{R} _2^{(0)}\cup \mathcal{R} _1\Gamma
_1\odot \mathcal{R} _2^{(1)}$\ and
$(\mathcal{R} _1 ')^{(0)}\odot \mathcal{R} _2 '\cup (\mathcal{R} _1
')^{(1)}\odot R_2 ' \Gamma _2$.\\
In order to prove that 
$(Z(\mathcal{R} _1)^{(0)}\overline{\otimes }Z(\mathcal{R} _2)^{(0)})'
\subseteq Z(\mathcal{R} _1\hat{\otimes } \mathcal{R} _2)'$, it suffices to show that $A\otimes B \in Z(\mathcal{R} _1\hat{\otimes } \mathcal{R} _2)'$ for each $A \in \mathcal{R} _1 \cup \mathcal{R} _1'\cup \{\Gamma _1\}''$,\ $B
\in \mathcal{R} _2 \cup \mathcal{R} _2'\cup \{\Gamma _2\}''$ because $Z(\mathcal{R} _1)^{(0)'}  \overline{\otimes } Z(\mathcal{R}
_2)^{(0)'}$ is generated by elements of this form.
It suffices to show that $\Gamma _1\otimes I$,\ $I \otimes \Gamma
_2$ are in $Z(\mathcal{R} _1\hat{\otimes } \mathcal{R} _2)'$ . In fact, since 
\begin{align}
(\mathcal{R} _1\Gamma _1\odot \mathcal{R} _2^{(1)})(\Gamma _1\otimes
I) &= \mathcal{R} _1\odot \mathcal{R} _2^{(1)},\\
((\mathcal{R} _1 ')^{(1)}\odot R_2 ' \Gamma _2)(I \otimes \Gamma _2)
&= (\mathcal{R} _1 ')^{(1)}\odot \mathcal{R} _2',
\end{align}
$\mathcal{R} _1\odot \mathcal{R} _2$ and $\mathcal{R} _1'\odot
\mathcal{R} _2'$ are included in $Z(\mathcal{R} _1\hat{\otimes } \mathcal{R} _2)'$ so $A\otimes B \in Z(\mathcal{R} _1\hat{\otimes } \mathcal{R}
_2)'$.\\
We shall show that $\Gamma _1\otimes I$ is included in $Z(\mathcal{R} _1\hat{\otimes } \mathcal{R}
_2)'$.\ Since $Z(\mathcal{R} _2)$ is balanced, there is a
self-adjoint\ unitary\ $U_2 \in \mathcal{R} _2^{(1)}\cap (\mathcal{R}
_2')^{(1)}$. Now 
\[\Gamma _1\otimes U_2 \in \mathcal{R} _1\Gamma _1\odot \mathcal{R} _2^{(1)}\]
\[I \otimes U_2 \in (\mathcal{R} _1 ')^{(0)}\odot \mathcal{R} _2 '\]
so
\[\Gamma _1\otimes I \in Z(\mathcal{R} _1\hat{\otimes } \mathcal{R} _2)'.\]\
Similarly $I \otimes \Gamma _2$ is in $Z(\mathcal{R} _1\hat{\otimes } \mathcal{R} _2)'$.$\Box $
\vspace{\baselineskip}

\begin{prop}\label{20}
Let $(\mathcal{R} _1,\ \text{Ad}_{\Gamma _1})$ and $(\mathcal{R} _2,\ \text{Ad}_{\Gamma
_2})$ be graded\ von\ Neumann\ algebras. If $(Z(\mathcal{R} _2),\ \text{Ad}_{\Gamma _2})$ is balanced, then 
$(\mathcal{R} _1\hat{\otimes }\mathcal{R} _2)^{(0)} $ is $*$-ismorphic to $\mathcal{R}
_1\overline{ \otimes} \mathcal{R} _2^{(0)}$.
\end{prop}
$<$proof$>$ Since $Z(\mathcal{R} _2)$ is balanced, there is a self-adjoint\ unitary $U_2 \in Z(\mathcal{R}
_2)^{(1)}$.\ Then 
$(\mathcal{R} _1\overline{\otimes } \mathcal{R} _2^{(0)} ,\ \text{Ad}_{\Gamma
_1\otimes U_2})$ becomes a graded\ von\ Neumann\ algebra. In fact, 
$\mathcal{R} _1\odot  \mathcal{R} _2^{(0)}$ is closed under $\text{Ad}_{\Gamma _1\otimes
U_2}$. We shall show that
\[(\mathcal{R} _1\overline{\otimes } \mathcal{R} _2^{(0)} )^{<1>} =
\overline{ \mathcal{R} _1^{(1)}\odot  \mathcal{R} _2^{(0)}},\]
where right-hand side is weak-operator closure of $\mathcal{R} _1^{(1)}\odot  \mathcal{R} _2^{(0)}$,\ and 
$(\mathcal{R} _1\overline{\otimes } \mathcal{R} _2^{(0)} )^{<0>}$, 
$(\mathcal{R} _1\overline{\otimes } \mathcal{R} _2^{(0)} )^{<1>}$ denote even, odd elements defined by grading operator $\text{Ad}_{\Gamma _1\otimes U_2}$ respectively. Since 
\[\mathcal{R} _1^{(1)}\odot  \mathcal{R} _2^{(0)} \subseteq
(\mathcal{R} _1\overline{\otimes } \mathcal{R} _2^{(0)} )^{<1>},\]it follows that 
\[\overline{\mathcal{R} _1^{(1)}\odot  \mathcal{R} _2^{(0)}} \subseteq
(\mathcal{R} _1\overline{\otimes } \mathcal{R} _2^{(0)} )^{<1>}.\]Assume that 
$A \in (\mathcal{R} _1\overline{\otimes } \mathcal{R} _2^{(0)} )^{<1>}$.\
 Let $\{\sum_i A_i^a\otimes B_i^a\}_a \subseteq \mathcal{R} _1\odot
\mathcal{R} _2^{(0)}$ be a net weak-operator converges to $A$.\ Since $U_2$ is a element of the center
\[ -\text{Ad}_{\Gamma _1\otimes U_2}(\sum_i A_i^a\otimes B_i^a)+\sum_i
A_i^a\otimes B_i^a \\
=2\sum_i (A_i^a)^{(1)} \otimes B_i^a \to 2A \ (\text{WOT}).\]
Thus $\sum_i (A_i^a)^{(1)}\otimes B_i^a \to A\ (\text{WOT})$ and $A \in
\overline{ \mathcal{R} _1^{(1)}\odot  \mathcal{R} _2^{(0)}}$.
Hence \[(\mathcal{R} _1\hat{\otimes }\mathcal{R} _2)^{<1>} = \overline{
\mathcal{R} _1^{(1)}\odot  \mathcal{R} _2^{(0)}}.\]
Similarly, we can show that $(\mathcal{R} _1\overline{\otimes }\mathcal{R} _2^{(0)})^{<0>} =
\mathcal{R} _1^{(0)}\overline{ \otimes} \mathcal{R} _2^{(0)}$.\\
Since $(\mathcal{R} _1\overline{\otimes } \mathcal{R} _2^{(0)} ,\
\text{Ad}_{\Gamma _1\otimes U_2})$ is a graded\ von\ Neumann\ algebra, 
$\mathcal{R} _1\overline{\otimes } \mathcal{R} _2^{(0)}$ is $*$-isomorphic to 
$\mathcal{R} _1^{(0)}\overline{ \otimes} \mathcal{R} _2^{(0)} \oplus
\overline{ \mathcal{R} _1^{(1)} \odot  \mathcal{R} _2^{(0)}}(\Gamma _1
\otimes U_2) $
 from \ref{n4}. $\mathcal{R} _1^{(0)}\overline{ \otimes} \mathcal{R} _2^{(0)} \oplus
\overline{ \mathcal{R} _1^{(1)} \odot  \mathcal{R} _2^{(0)}}(\Gamma _1
\otimes U_2) $ is generated by $\mathcal{R} _1^{(0)}\odot  \mathcal{R} _2^{(0)}$ and 
$\mathcal{R} _1^{(1)}\Gamma_1\odot
\mathcal{R}_2^{(1)}$ so that it equals to $(\mathcal{R}_1\hat{\otimes}\mathcal{R}_2)^{(0)}$.$\Box
$
\vspace{\baselineskip}

\begin{prop}\label{21}
Let $(\mathcal{R}_1,\ \text{Ad}_{\Gamma_1})$ and $ (\mathcal{R}_2,\
\text{Ad}_{\Gamma_2})$ be central\ graded\ von\ Neumann\
algebras and not factors.\ 
If $\mathcal{R}_1\overline{\otimes}\mathcal{R}_2^{(0)} $ is of type\
I\hspace{-1.2pt}I$_\infty$, $\mathcal{R}_1\hat{\otimes}\mathcal{R}_2$ is also of I\hspace{-1.2pt}I$_\infty$.
\end{prop}
$<$proof$>$ Since $Z(\mathcal{R} _2)$ is balanced, $Z((\mathcal{R}
_1\hat{\otimes } \mathcal{R} _2)^{(0)})$ is $*$-isomorphic to 
$Z(\mathcal{R} _1\overline{\otimes }\mathcal{R} _2^{(0)})$ from \ref{20}.
Thus $Z((\mathcal{R} _1\hat{\otimes } \mathcal{R} _2)^{(0)})$ is finite-dimensional from \ref{17}.\ Furthermore, 
$Z(\mathcal{R} _2)$ is balanced so the type of 
$\mathcal{R} _1\overline{\otimes }\mathcal{R}
_2^{(0)}$ is equal to the type of $(\mathcal{R} _1\hat{\otimes } \mathcal{R}
_2)^{(0)}$ from \ref{20}.\ Since both
$(\mathcal{R}_1,\ \text{Ad}_{\Gamma_1})$ and $(\mathcal{R}_2,\
\text{Ad}_{\Gamma_2})$ are central, 
$\mathcal{R}_1\hat{\otimes}\mathcal{R}_2$ is a factor from \ref{18}. Thus  $\mathcal{R}_1\hat{\otimes}\mathcal{R}_2$ is of type\ $\mu$ ($\mu$
=I,\ I\hspace{-1.2pt}I,\ I\hspace{-1.2pt}I\hspace{-1.2pt}I).\ Since $\mathcal{R}_1\hat{\otimes}\mathcal{R}_2$ is 
balanced, the type of $\mathcal{R}_1\hat{\otimes}\mathcal{R}_2$ is equal to the type of $(\mathcal{R} _1\hat{\otimes }
\mathcal{R} _2)^{(0)}$ from \cite{bobo} Lem A.1.\ $I$ is an infinite relative to   
$(\mathcal{R} _1\hat{\otimes } \mathcal{R}
_2)^{(0)}$ so $I$ is infinite relative to $\mathcal{R}_1\hat{\otimes}\mathcal{R}_2$. Thus 
$\mathcal{R}_1\hat{\otimes}\mathcal{R}_2$ is of type\ I\hspace{-1.2pt}I$_{\infty}$.$\Box $
\vspace{\baselineskip}

We can show the following proposition by using topological characterization of finite projections.  

\begin{prop}\label{35}
Let $(\mathcal{R} ,\ \text{Ad}_\Gamma )$ be a graded\ von\ Neumann\ algebra which has an odd unitary\ $U$.\ Let $E \in \mathcal{R}
^{(0)}$ be a projection which is finite relative to $\mathcal{R} ^{(0)}$. Then $E$ is finite relative to $\mathcal{R} $.
\end{prop}
$<$proof$>$ It suffices to show that the map $A \mapsto EA^*$,\ $\mathcal{R}
\rightarrow \mathcal{R} $ is strong-operator continuous on the unit ball from \cite{krkr} 11.2.23.\ Let 
$\{A_a\}_a$ be a net such that $A_a,\  A\in  (\mathcal{R} )_1$, $A_a\to A$ (SOT).\ Since $\text{Ad}_{\Gamma }$ is strong-operator continuous, it follows that $A_a^{(0)} \to A^{(0)}$,\ $A_a^{(1)} \to A^{(1)}$ (SOT) and they are element of $(\mathcal{R} )_1$.\ Since $E$ is finite relative to $\mathcal{R}^{(0)}$, 
$E(A_a^{(0)})^* \to E(A^{(0)})^*$,\
$E(A_a^{(1)})^*U^* \to E(A^{(1)})^*U^*$ (SOT) and 
$E(A_a^{(1)})^* \to E(A^{(1)})^*$ (SOT) from \cite{krkr} 11.2.23. Thus $EA_a^* \to
EA$ (SOT) and $E$ is finite relative to $\mathcal{R} $.
$\Box $
\vspace{\baselineskip}

\begin{prop}\label{23}
Let $(\mathcal{R},\ \text{Ad}_\Gamma)$ be a graded\ von\ Neumann\
algebra such that $(Z(\mathcal{R}) ,\ \text{Ad}_\Gamma )$ is balanced.\ If $\mathcal{R}$ is of type\ I\hspace{-1.2pt}I$_1$,\ I\hspace{-1.2pt}I$_\infty$, then $\mathcal{R}^{(0)}$ is also of \\type\ I\hspace{-1.2pt}I$_1$,\ I\hspace{-1.2pt}I$_\infty$ respectively.
\end{prop}
$<$proof$>$ In either case, $\mathcal{R} ^{(0)}$ has no non-zero abelian\ projection from \ref{19} because $(Z(\mathcal{R}) ,\ \text{Ad}_\Gamma)$ is balanced.\ Assume that $\mathcal{R} $ is of type\ I\hspace{-1.2pt}I$_1$.\ We note that $\mathcal{R}^{(0)}$ has no abelian\ projection and finite so it is of type\ I\hspace{-1.2pt}I$_1$.\\
\ Assume that $\mathcal{R} $ is of type\ I\hspace{-1.2pt}I$_\infty$.\ Since $\mathcal{R} $ is of type\
I\hspace{-1.2pt}I, there is a finite projection $F \in \mathcal{R} $ in $\mathcal{R} $
such that its central\ support $C_F$ in $\mathcal{R} $ is $ I$. We note that
\[\text{Ad}_{\Gamma}(F\lor \text{Ad}_{\Gamma}(F) )=\text{Ad}_{\Gamma}(F)\lor
\text{Ad}_{\Gamma}(\text{Ad}_{\Gamma}(F))=F\lor \text{Ad}_{\Gamma}(F)\]and 
$F\lor \text{Ad}_{\Gamma}(F) $ is even.\ Since $F$ is finite relative to $\mathcal{R}
$, $\text{Ad}_{\Gamma}(F)=\Gamma F \Gamma $ is finite relative to $\mathcal{R}
$.\ Thus $F\lor \text{Ad}_{\Gamma}(F) $ is finite relative to $\mathcal{R} $ from \cite{krkr} 6.3.8. THEOREM. Assume that there is a non-zero projection $Q \in Z(\mathcal{R} ^{(0)})$ such that $\mathcal{R}
^{(0)}Q$ is of type\ I\hspace{-1.2pt}I\hspace{-1.2pt}I.\ Since $Q,\ F$ have central\ carriers in $\mathcal{R} $ which is $C_Q C_F\neq 0$, they have an equivalent non-zero subprojection in $\mathcal{R} $ from \cite{krkr} 6.1.8. PROPOSITION.\ Let $G$ be a non-zero subprojection of $Q$ which is equivalent to a subprojection of $F$. $G\vee \Gamma
G\Gamma$ is even and finite relative to $\mathcal{R} $ from the fact we have shown.\ Since $G\leq Q$ and $Q$ is even, $G\vee \Gamma G\Gamma\leq Q$ and $G\vee \Gamma
G\Gamma$ is also finite relative to $\mathcal{R} ^{(0)}Q$, contradicting the fact that $\mathcal{R}^{(0)}Q$ is of type I\hspace{-1.2pt}I\hspace{-1.2pt}I.
Thus $\mathcal{R} ^{(0)}$ has no central portion of type\ I\hspace{-1.2pt}I\hspace{-1.2pt}I.\\
Furthermore, abelian\ projections in $\mathcal{R} ^{(0)}$ are also abelian in $\mathcal{R}
$, so $\mathcal{R} ^{(0)}$ has no central portion of type\ I.\\
So far, we have proved $\mathcal{R} ^{(0)}$ is of type\ I\hspace{-1.2pt}I.\ Assume that\ there is a
central projection\ $P_{c_1}\in Z(\mathcal{R} ^{(0)})$ such that
$\mathcal{R}^{(0)}P_{c_1}$ is of type\ I\hspace{-1.2pt}I$_1$.We note that $P_{c_1}\in Z(\mathcal{R})$ from \ref{17} because $(Z(\mathcal{R}),\ \text{Ad}_\Gamma )$ is balanced.\ Since $(\mathcal{R}P_{c_1} ,\
\text{Ad}_{P_{c_1}\Gamma} )$ has an odd unitary\ $UP_{c_1}$, 
$P_{c_1}$ is finite in $\mathcal{R}P_{c_1}$ from \ref{35}.\ We note that $\mathcal{R} $ is of type\ I\hspace{-1.2pt}I$_\infty$, contradicting the fact that $\mathcal{R}P_{c_1}$ is finite.
Thus $\mathcal{R} ^{(0)}$ is of type\ I\hspace{-1.2pt}I$_\infty$.$\Box $
\vspace{\baselineskip}

\begin{prop}\label{21-1}
Let $(\mathcal{R} _1,\ \text{Ad}_{\Gamma_1})$ be a central graded von Neumann algebra such that  $\mathcal{R}_1$ is of type\ I\hspace{-1.2pt}I$_1$ and not a factor,\ $(\mathcal{R} _2,\ \text{Ad}_{\Gamma_2})$ be a central graded von Neumann algebra such that $\mathcal{R} _2$ is of type\ I$_n$($n <\infty$) and not a factor. Then $\mathcal{R}_1 \hat{\otimes } \mathcal{R} _2$ is of type\ I\hspace{-1.2pt}I$_1$.
\end{prop}
$<$proof$>$ Since $\mathcal{R}  _2$ is central and of type\
I$_n$, $\mathcal{R} _2^{(0)}$ is of type\ I$_n$\ from \ref{11} and $\mathcal{R} _1 \overline{\otimes } \mathcal{R} _2^{(0)}$ is of type\
I\hspace{-1.2pt}I$_1$.\ $(Z(\mathcal{R}_2),\
\text{Ad}_{\Gamma_2})$ is balanced, so 
$(\mathcal{R} _1 \hat{\otimes } \mathcal{R} _2)^{(0)}$ is also of type\
I\hspace{-1.2pt}I$_1$ from \ref{20}.\ Since both $(Z(\mathcal{R}_1),\ \text{Ad}_{\Gamma_1})$ and $(Z(\mathcal{R}_2),\ \text{Ad}_{\Gamma_2})$ are balanced, it follows that $Z(\mathcal{R} _1 \hat{\otimes
} \mathcal{R} _2)=\mathbb{C} I$ from \ref{18}.\ The right-hand side of the isomorphism  $Z((\mathcal{R} _1 \hat{\otimes } \mathcal{R} _2)^{(0)})\cong
Z(\mathcal{R}_1\overline{\otimes}\mathcal{R}_2^{(0)})$ is finite-dimensional, so the left-hand side is also finite-dimensional.\ We note that $\mathcal{R} _1 \hat{\otimes } \mathcal{R}
_2$ is balanced because $(\mathcal{R} _1,\ \text{Ad}_{\Gamma_1})$ is central and not a factor. Thus the type of $(\mathcal{R} _1 \hat{\otimes } \mathcal{R}
_2)^{(0)}$ is equal to the type of $\mathcal{R} _1 \hat{\otimes } \mathcal{R} _2$ from 
\cite{bobo} Lemma.A.1.\ In addition, $\mathcal{R} _1 \hat{\otimes }
\mathcal{R} _2$ is finite from \ref{35} because it is balanced.
Thus $\mathcal{R} _1 \hat{\otimes }
\mathcal{R} _2$ is of type\ I\hspace{-1.2pt}I$_1$.$\Box $
\vspace{\baselineskip}

\begin{prop}\label{49}
Let $(\mathcal{R}_1,\ \text{Ad}_{\Gamma_1})$ and $(\mathcal{R}_2,\
\text{Ad}_{\Gamma_2})$ be central\ graded\ von\ Neumann\
algebras which are not factors.\
Then multiplication table of types of graded\ tensor\ product $\mathcal{R}_1 \hat{\otimes } \mathcal{R} _2$ coincides with  that of normal tensor product $\mathcal{R}_1 \overline{\otimes } \mathcal{R} _2$ except the case that both $\mathcal{R}_1$,\
$\mathcal{R}_2$ are of type\ I.\ In the case
$\mathcal{R} _1$ is of type\ I$_m$($m$ is a cardinal number) and $\mathcal{R} _2$ is of type\
I$_n$($n$ is a cardinal number), their graded tensor product is of type\ I$_{2mn}$.
\end{prop}
$<$proof$>$ In the case at least one of $\mathcal{R}_1$ and $\mathcal{R}_2$ is of type\ I\hspace{-1.2pt}I\hspace{-1.2pt}I,\ and in the case both $\mathcal{R}_1$ and $\mathcal{R}_2$ is of type\ I, we have already proved coincidence in section 3.1,\ Proposition \ref{12} respectively.\ So we shall consider other cases. If $\mathcal{R}_1$ is of type\ I\hspace{-1.2pt}I$_\infty$, the type of $\mathcal{R}_1 \hat{\otimes } \mathcal{R} _2$ is also of type\
I\hspace{-1.2pt}I$_\infty$ from \ref{20}. In the case $\mathcal{R}_2$ is of type\ I\hspace{-1.2pt}I$_\infty$, the type of $\mathcal{R}_1 \hat{\otimes } \mathcal{R} _2$ is also of type\ I\hspace{-1.2pt}I$_\infty$ from \ref{48} similarly.\ $\mathcal{R} _1 \hat{\otimes } \mathcal{R} _2  $ is a factor because both center of $\mathcal{R} _1$ and center of $\mathcal{R} _2$ are 
balanced.\ In addition, since $Z(\mathcal{R}_2)$ is balanced, $(\mathcal{R}_1\hat{\otimes}\mathcal{R}_2)^{(0)}
$ is of type I or I\hspace{-1.2pt}I and the center of $(\mathcal{R}_1\hat{\otimes}\mathcal{R}_2)^{(0)}
$ is finite-dimensional from the $*$-isomorphism of \ref{20}.
Thus we can apply \cite{bobo} Lemma.A.1 to $\mathcal{R}_1\hat{\otimes}\mathcal{R}_2$ from \ref{18}, and the type (=I, I\hspace{-1.2pt}I) of $\mathcal{R} _1\hat{\otimes } \mathcal{R} _2  $ is equal to the type(=I, I\hspace{-1.2pt}I) of 
$(\mathcal{R} _1 \hat{\otimes } \mathcal{R} _2)^{(0)}  $.\ 
$\mathcal{R}_1\hat{\otimes}\mathcal{R}_2 $ is finite if and only if $(\mathcal{R} _1 \hat{\otimes } \mathcal{R} _2)^{(0)}  $ is finite from \ref{35} because $\mathcal{R}_1\hat{\otimes}\mathcal{R}_2 $ is balanced.\ Since $(Z(\mathcal{R}_2),\
\text{Ad}_{\Gamma_2})$ is balanced, the type of $(\mathcal{R} _1 \hat{\otimes } \mathcal{R} _2)^{(0)}  $ is equal to the type of $\mathcal{R} _1 \overline{\otimes } (\mathcal{R} _2)^{(0)} $ from \ref{20}. Thus we shall check the type of $\mathcal{R} _1\overline{ \otimes} \mathcal{R} _2^{(0)}$. We note that the type of $\mathcal{R} _2^{(0)}$ is equal to the type of $\mathcal{R} _2$  from \ref{11},\ \ref{23} and the
type of $\mathcal{R} _1\overline{ \otimes} \mathcal{R} _2^{(0)}$ is equal to $\mathcal{R} _1\overline{ \otimes} \mathcal{R} _2$.$\Box $

\vspace{\baselineskip}

\subsection{The case $(\mathcal{R}_2,\
\text{Ad}_{\Gamma_2})$ is central, $\mathcal{R}_1$ is a factor, and
$\mathcal{R}_2$ is not a factor}

\begin{prop}\label{2626}
Every $*$-automorphism of $\mathfrak{B}  (\mathcal{H} )$ is inner. That is, let $\varphi : \mathfrak{B} ( \mathcal{H} ) \twoheadrightarrow
\mathfrak{B} ( \mathcal{H} ) $ be a $*$-automorphism, then there is a unitary\ $U\in \mathfrak{B} (\mathcal{H} )$ such that 
\[\varphi (A) = UAU^*\]
for each $A \in \mathfrak{B}(\mathcal{H} )$.
\end{prop}
$<$proof$>$ \cite{inin} Example 2.6.26.$\Box $
\vspace{\baselineskip}

\begin{prop}\label{8}
If $\varphi $ is a $*$-automorphism\ of $\mathfrak{B} (\mathcal{H} )$ such that $\varphi
\circ \varphi = \text{id} $, there is a self-adjoint\ unitary
$U \in \mathfrak{B} (\mathcal{H} )$ such that\[\varphi (A) = UAU \]
for each $A \in \mathfrak{B}(\mathcal{H} )$.
\end{prop}
$<$proof$>$ There is a unitary\ $U \in \mathfrak{B} (\mathcal{H})$ such that
\[\varphi (A) = UAU^*\]
for each $A \in \mathfrak{B} (\mathcal{H} )$ from \ref{2626}.\ Since $\varphi \circ \varphi = \text{id} $,
\[U^2A{U^*}^2 = A,\:U^2A = AU^2.\]
for each $A \in \mathfrak{B} (\mathcal{H} )$. Let $\mathcal{C} $ be a center of $\mathfrak{B}(\mathcal{H} )$ then $U^2 \in \mathcal{C} $.
Suppose that $U^2 = \alpha I,\ \alpha \in \mathbb{C} $. We note that 
the norm of $U^2$ is $1$ and $|\alpha | = 1$, so that there is a $\theta\in \mathbb{R}$ such that $\alpha  = e^{i\theta} ,\ \theta
\in \mathbb{R} $.\ Assume that $U' = e^{-i\theta /2}U$ then it follows that
\[U'^* = e^{i\theta /2} U^*= e^{i\theta /2}e^{-i\theta }U= U'.\]
Consequently $U'$ is self-adjoint.\ $U'$ is also a unitary. In addition, we note that 
\[U'AU'^* = UAU^*= \varphi (A)\] for each $A \in \mathfrak{B}
(\mathcal{H} )$. Thus $U'$ satisfies required conditions.$\Box$
\vspace{\baselineskip}

\begin{prop}\label{9}
Let $(\mathcal{R} ,\ \text{Ad}_\Gamma )$ be a spatially\ graded\ von\ Neumann\
algebra\ such that $\mathcal{R} $ is a type I\ factor. Then there are a Hilbert space $\mathcal{H}$ and a self-adjoint unitary $\Gamma$ such that
$(\mathcal{R} ,\ \text{Ad}_\Gamma )$ is graded $*$-isomorphic to $(\mathfrak{B} (\mathcal{H} ),\ \text{Ad}_U)$.
\end{prop}
$<$proof$>$ Since $\mathcal{R} $ is a type I\ factor, there is a Hilbert space $\mathcal{H}$ such that $\mathcal{R} $ is $*$-isomophic to $\mathfrak{B} (\mathcal{H})$  from \cite{krkr}\ 6.6.1.THEOREM.\ Let 
$\varphi :\mathcal{R} \twoheadrightarrow \mathfrak{B} (\mathcal{H})$ be a $*$-isomorphism. We can define a grading on $\mathfrak{B} (\mathcal{H} )$ by $\varphi \circ \text{Ad}_\Gamma \circ\varphi ^{-1}$. Thus there is a self-adjoint\ unitary\ $U \in \mathfrak{B} (\mathcal{H} )$ such that 
\[ (\varphi \circ \text{Ad}_\Gamma \circ \varphi ^{-1} )(A)= \text{Ad}_U (A)\] for all $A \in\mathfrak{B} (\mathcal{H} )$ from \ref{8}. We shall show that the map\ $\varphi : (\mathcal{R} ,\ \text{Ad}_\Gamma )\twoheadrightarrow
(\mathfrak{B} (\mathcal{H} ),\ \text{Ad}_U)$ is a graded\ $*$-isomorphism.\ It suffices to show that $\text{Ad}_\Gamma$ preserve the grading. Now,
\[(\text{Ad}_U \circ \varphi) (A) = (\varphi \circ \text{Ad}_\Gamma )(A)= \varphi (A)\]
for all $A \in\mathcal{R} ^{(0)}$ so $\varphi (A) \in \mathfrak{B} (\mathcal{H} )^{(0)}$.\ It will be similarly shown that $
\varphi(\mathcal{R} ^{(1)}) \subseteq \mathfrak{B} (\mathcal{H} )^{(1)}$.$\Box$
\vspace{\baselineskip}

\begin{prop}\label{10}
Let $(\mathcal{R} _1,\ \text{Ad}_{\Gamma _1}) $ be a graded\ von\ Neumann\
algebra such that $\mathcal{R} _1$ is a type I\ factor and 
\ $(\mathcal{R} _2,\ \text{Ad}_{\Gamma _2}) $ be a balanced\ graded\ von\ Neumann\ algebra. Then $\mathcal{R} _1 \hat{\otimes  } \mathcal{R} _2$ is $*$-isomorphic to $\mathcal{R} _1\overline{ \otimes} \mathcal{R} _2 $.
\end{prop}
$<$proof$>$ By \ref{9}, there are a Hilbert space $\mathcal{H}$ and a self-adjoint unitary $U\in \mathfrak{B} (\mathcal{H} )$ such that $(\mathfrak{B} (\mathcal{H} ),\ \text{Ad}_U)$,\ is graded $*$-isomorphic to 
$(\mathcal{R} _1,\ \text{Ad}_{\Gamma _1}) $.\ Since 
$U \in \mathfrak{B} (\mathcal{H} )$, $\mathfrak{B} (\mathcal{H} )
\hat{\otimes } \mathcal{R} _2$ is $*$-isomorphic to 
$ \mathfrak{B} (\mathcal{H} ) \overline{ \otimes} \mathcal{R} _2$
. In fact,\ $AU^{\partial B} \otimes B \in \mathfrak{B} (\mathcal{H} ) \overline{\otimes} \mathcal{R} _2 $ for homogeneous elementary\ tensors so
\[ \mathfrak{B} (\mathcal{H})\hat{\otimes}\mathcal{R} _2
\subseteq \mathfrak{B} (\mathcal{H} ) \overline{ \otimes} \mathcal{R}_2.\]
Conversely, assume that $A \in \mathfrak{B} (\mathcal{H} )$,\ $B \in \mathcal{R} _2$ then 
\[A\otimes B = A^{(0)}\otimes B^{(0)}+A^{(1)}\otimes B^{(0)}
+A^{(0)}UU\otimes B^{(1)}+A^{(1)}UU\otimes B^{(1)}\]
and $A\otimes B \in  \mathfrak{B} (\mathcal{H} ) \hat{\otimes }
\mathcal{R} _2$. Thus
\[ \mathfrak{B} (\mathcal{H} ) \hat{\otimes } \mathcal{R} _2=
\mathfrak{B} (\mathcal{H} ) \overline{ \otimes} \mathcal{R} _2 .\]
\ Since $(\mathcal{R} _2,\ \text{Ad}_{\Gamma _2})
$ is balanced, 
\[\mathcal{R} _1 \hat{\otimes  } \mathcal{R} _2 \cong \mathfrak{B}
(\mathcal{H} ) \hat{\otimes } \mathcal{R} _2
\cong  \mathfrak{B} (\mathcal{H} ) \overline{ \otimes} \mathcal{R}
_2\cong  \mathcal{R} _1 \overline{ \otimes} \mathcal{R} _2\]
from \cite{bobo} Lemma A.5.$\Box$
\vspace{\baselineskip}

\begin{prop}
Let $(\mathcal{R}_1,\ \text{Ad}_{\Gamma_1})$ be a graded von Neumann algebra such that $\mathcal{R}_1$ is a type\ I\ factor, and $(\mathcal{R}_2,\ \text{Ad}_{\Gamma_2})$ be a central graded von Neumann algebra such that $\mathcal{R}_2$ is not a factor. Then the type of $\mathcal{R} _1 \hat{\otimes } \mathcal{R} _2  $ is equal to  the type of $\mathcal{R} _1 \overline{ \otimes} \mathcal{R} _2 $.
\end{prop}
$<$proof$>$ Since $(\mathcal{R}_2,\ \text{Ad}_{\Gamma_2})$ is central and $\mathcal{R}_1$ is a type\ I\ factor, $\mathcal{R} _1 \hat{\otimes  } \mathcal{R} _2$ is $*$-isomorphic to $\mathcal{R} _1\overline{ \otimes} \mathcal{R} _2 $ from \ref{10}. We note that the types of von Neumann algebras are preserved by $*$-isomorphisms.$\Box$

\vspace{\baselineskip}

\begin{prop}\label{34}
Let $\mathcal{R} _1$ and $\mathcal{R} _2$ be von\ Neumann\ algebras on a Hilbert space $\mathcal{H} $. If
$\mathcal{R} _1\subseteq \mathcal{R} _2$, $\mathcal{R} _1$ is of type\
I\hspace{-1.2pt}I$_1$ and $\mathcal{R} _2$ is finite, then 
$\mathcal{R} _2$ is also of type\ I\hspace{-1.2pt}I$_1$.
\end{prop}
$<$proof$>$ Assume that $P\in Z(\mathcal{R} _2)$ is a non-zero projection such that $\mathcal{R}_2P$ is of type\ I$_n$.\ Since  $\mathcal{R} _1$ is of type\ I\hspace{-1.2pt}I$_1$, there are projections $E_k \in \mathcal{R} _1$ ($k=1,2,\dots,n+1$) such that $I=E_1+...+E_{n+1}$ and $E_i\sim E_j$ in $\mathcal{R}_1$ from \cite{krkr} 6.5.6. LEMMA.\ We note that 
$E_i\sim E_j$ holds true in $\mathcal{R} _2$.\ Since a central\ carrier preserves unions from \cite{krkr} 5.5.3.\ PROPOSITION, central carrier $C_{E_i}$ is equal to $I$ in $\mathcal{R} _2$, for all $i$. We note that $PE_i\neq0$,\ $PE_i\sim
PE_j$.\ Since $\mathcal{R} _2P$ is of type\ I$_n$, there is a non-zero abelian\ projection\ $F_1$ such that $F_1\leq PE_1$.\ Let $V_j$ be a partial isometry with initial projection $E_1$ and final projection $E_{j}$, then $V_jF_1$ is a partial isometry with initial projection $F_1$ and final projection $E_j'$, where $E_j'$  is a subprojection of $E_j$. Thus there are $n+1$ equivalent abelian\ projections in $\mathcal{R}_2P$, contradicting the fact that $\mathcal{R} _2P$ is of type\ I$_n$ and \cite{krkr} 6.5.2. THEOREM(Type decomposition).$\Box$
\vspace{\baselineskip}

\begin{prop}\label{28}
Let $(\mathcal{R} _1,\ \text{Ad}_{\Gamma _1})$ and $(\mathcal{R} _2,\ \text{Ad}_{\Gamma
_2})$ be graded\ von\ Neumann\ algebras\ such that 
$\mathcal{R} _1$ is a factor,\ $(\mathcal{R} _2,\ \text{Ad}_{\Gamma
_2})$ is central and $\mathcal{R} _2$ is not a factor. Furthermore, if at least one of $\mathcal{R}_1$,\ $\mathcal{R} _2$ is of type\ I\hspace{-1.2pt}I$_1$\ and the other is finite,\ then $\mathcal{R}
_1 \hat{\otimes } \mathcal{R} _2  $ is of type\ I\hspace{-1.2pt}I$_1$.
\end{prop}
$<$proof$>$ Let $U_2 \in Z(\mathcal{R} _2)$ be a self\ adjoint\ unitary.\ We note that 
$\Gamma _1 \otimes U_2$ is an odd unitary in $(\mathcal{R} _1 \hat{\otimes } \mathcal{R} _2,\ \text{Ad}_{I \otimes
\Gamma_2})$.\ With the notation of \ref{32}
\[(\mathcal{R} _1 \hat{\otimes } \mathcal{R} _2)^{<0>}=\mathcal{R}
_1\overline{\otimes }(\mathcal{R} _2)^{(0)} \]
from \ref{32}.\ Since $(\mathcal{R} _2,\ \text{Ad}_{\Gamma
_2})$ is central and  $\mathcal{R}_2$ is not a factor, $\mathcal{R}
_1\overline{\otimes }(\mathcal{R} _2)^{(0)}$ is of type I\hspace{-1.2pt}I$_1$ from \ref{23}.\ $\mathcal{R} _1 \hat{\otimes } \mathcal{R}_2$ is finite from \ref{35} because it has an odd unitary.\ Thus $\mathcal{R} _1 \hat{\otimes } \mathcal{R} _2$ is of type\ I\hspace{-1.2pt}I$_1$ from \ref{34}.$\Box $
\vspace{\baselineskip}

\begin{prop}\label{26}
Let $(\mathcal{R},\ \text{Ad}_\Gamma)$ be a graded\ von\ Neumann\ algebra such that $\mathcal{R}^{(0)}$ is of type\ I\hspace{-1.2pt}I\hspace{-1.2pt}I. Then $\mathcal{R}$ is of type\ I\hspace{-1.2pt}I\hspace{-1.2pt}I.
\end{prop}
$<$proof$>$ If there is a non-zero projection $F\in \mathcal{R} $ which is finite relative to $\mathcal{R} $,\ then $F\lor \text{Ad}_{\Gamma}(F) $ is even.\ Since $F$ is finite relative to $\mathcal{R} $, $\Gamma
F \Gamma $ is also finite relative to $\mathcal{R} $. Thus $F\lor  \Gamma F\Gamma \in\mathcal{R}^{(0)}$ is finite relative to $\mathcal{R} $ so that it is finite relative to $\mathcal{R}^{(0)}$ from \cite{krkr} 6.3.8\ THEOREM. But this contradicts the fact that $\mathcal{R}^{(0)}$ is of type\ I\hspace{-1.2pt}I\hspace{-1.2pt}I.$\Box$
\vspace{\baselineskip}

\begin{prop}\label{27}
Let $(\mathcal{R} ,\ \text{Ad}_\Gamma )$ be a graded\ von\ Neumann\
algebra which has an odd unitary. If $\mathcal{R} $ is of type\ I\hspace{-1.2pt}I\hspace{-1.2pt}I, then $\mathcal{R} ^{(0)}$ is also of type\ I\hspace{-1.2pt}I\hspace{-1.2pt}I.
\end{prop}
$<$proof$>$ Assume that there is a non-zero projection $E \in
\mathcal{R} ^{(0)}$ which is finite relative to $\mathcal{R} ^{(0)}$.\ We note that $\mathcal{R} $ has an odd unitary then $E$ is finite relative to $\mathcal{R}
$ from \ref{35}, contradicting the fact that $\mathcal{R}$ is of type\ I\hspace{-1.2pt}I\hspace{-1.2pt}I. Thus $\mathcal{R} ^{(0)}$ has no finite\ projection and $\mathcal{R} ^{(0)}$ is of type\ I\hspace{-1.2pt}I\hspace{-1.2pt}I.$\Box $
\vspace{\baselineskip}

\begin{prop}\label{36}
Let $(\mathcal{R} _1,\ \text{Ad}_{\Gamma _1})$ and $(\mathcal{R} _2,\ \text{Ad}_{\Gamma
_2})$ be graded\ von\ Neumann\ algebras, such that at least one of 
$\mathcal{R} _1$,\ $(\mathcal{R} _2)^{(0)}$ is properly\ infinite. Then 
$\mathcal{R} _1 \hat{\otimes } \mathcal{R} _2$ is also properly\ infinite.
\end{prop}
$<$proof$>$ Assume that $\mathcal{R} _1$ is properly\ infinite. There is a projection\ $E\in \mathcal{R} _1$ such that $I\sim  E\sim
I-E$ from \cite{krkr} 6.3.3.\ LEMMA.\  Let $V\in\mathcal{R}_1$ be a partial\ isometry\ with initial projection $E$ and final projection $I-E$, and $Q\in Z(\mathcal{R}_1 \hat{\otimes} \mathcal{R}_2)$ be a projection, then $Q(V\otimes I)$ is a partial isometry with initial projection $Q(E\otimes I)$ and final projection $Q[(I-E)\otimes I]$. Similarly it follows that $Q\sim Q(E\otimes I)$. Thus 
\[Q[(I-E)\otimes I]\sim Q(E\otimes I)\sim Q=Q[(I-E)\otimes I]+
Q(E\otimes I)\]and $I$ is properly\ infinite.\ From this same discussion, we can prove that, if $(\mathcal{R} _2)^{(0)}$ is properly\ infinite, $\mathcal{R} _1 \hat{\otimes } \mathcal{R} _2$ is also properly infinite.$\Box $
\vspace{\baselineskip}

\begin{prop}\label{37}
Let $(\mathcal{R} _1,\ \text{Ad}_{\Gamma _1})$ and $(\mathcal{R} _2,\ \text{Ad}_{\Gamma_2})$ be graded\ von\ Neumann\ algebras such that $\mathcal{R} _1$ is a factor,\ $(\mathcal{R} _2,\ \text{Ad}_{\Gamma_2})$ is central and $\mathcal{R} _2$ is not a factor.\ If either $\mathcal{R} _1$ is of type\ I\hspace{-1.2pt}I$_1$, $\mathcal{R} _2$ is infinite and not of type\ I\hspace{-1.2pt}I\hspace{-1.2pt}I \ or $\mathcal{R} _1$ is of  type I\hspace{-1.2pt}I$_\infty$ and $\mathcal{R}_2$ is of type $\mu$ ($\mu=$I, I\hspace{-1.2pt}I),\ then $\mathcal{R} _1 \hat{\otimes } \mathcal{R}_2$ is of type\ I\hspace{-1.2pt}I$_\infty$.
\end{prop}
$<$proof$>$ From \ref{17}, $(\mathcal{R}_2)^{(0)}$ is a factor. Thus $(\mathcal{R}_2)^{(0)}$ is of type I$_n$($n$ is a cardinal number) or I\hspace{-1.2pt}I$_1$, or I\hspace{-1.2pt}I$_\infty$ from \ref{26}.
So there are finite projections $E \in \mathcal{R} _1$,\ $F \in (\mathcal{R} _2)^{(0)}$ such that their central\ carriers are equal to $I$ in $\mathcal{R} _1$,\ $(\mathcal{R} _2)^{(0)}$ respectively.
$E\otimes F$ is a projection which is finite relative to $\mathcal{R} _1\overline{\otimes }(\mathcal{R}_2)^{(0)}$ and its central carrier in $\mathcal{R}_1\overline{\otimes}(\mathcal{R}_2)^{(0)}$ is $I$.\ Since $\mathcal{R} _1 \hat{\otimes }\mathcal{R} _2$ is balanced, $E\otimes F$ is also finite relative to $\mathcal{R} _1 \hat{\otimes } \mathcal{R} _2$ from \ref{35}. Moreover, since
\[\mathcal{R} _1 \hat{\otimes } \mathcal{R} _2 \supseteq \mathcal{R}
_1\overline{\otimes }(\mathcal{R} _2)^{(0)}, \]
 its central carrier is also $I$ in $\mathcal{R} _1 \hat{\otimes }
\mathcal{R} _2$ from \cite{krkr} 5.5.2. PROPOSITION. The right-hand side of the isomorphism
\[(E\otimes F)(\mathcal{R} _1\overline{\otimes }(\mathcal{R}
_2)^{(0)})(E\otimes F) \cong
E\mathcal{R} _1 E\overline{\otimes }F(\mathcal{R} _2)^{(0)}F\]
is of type\ I\hspace{-1.2pt}I$_1$ and $(E\otimes F)(\mathcal{R} _1 \hat{\otimes }
\mathcal{R} _2)(E\otimes F)$ is finite, so that 
 this is of type\ I\hspace{-1.2pt}I$_1$ from \ref{34}. Thus if $\mathcal{R} _1 \hat{\otimes } \mathcal{R}
_2$ has a non-zero abelian\ projection, it has a non-zero subprojection equivalent to a subprojection of $E\otimes
F$ from \cite{krkr} 6.1.8. PROPOSITION. It follows that 
$\mathcal{R} _1 \hat{\otimes } \mathcal{R} _2$ is of type\ I\hspace{-1.2pt}I.\ Since at least one of 
$\mathcal{R} _1$,\ $(\mathcal{R} _2)^{(0)}$ is properly\
infinite, $\mathcal{R} _1 \hat{\otimes } \mathcal{R}
_2$ is also properly infinite and is of type\ I\hspace{-1.2pt}I$_\infty$ from \ref{36}.$\Box $
\vspace{\baselineskip}

\begin{prop}\label{50}
Let $(\mathcal{R}_1,\ \text{Ad}_{\Gamma_1}),\ (\mathcal{R}_2,\
\text{Ad}_{\Gamma_2})$ be graded\ von\ Neumann\ algebras. If $(\mathcal{R}_2,\ \text{Ad}_{\Gamma_2})$ is central, $\mathcal{R}_1$ is a factor and $\mathcal{R}_2$ is not a factor, then multiplication table of graded\ tensor\
product is the same as that of normal tensor\ product.
\end{prop}
$<$proof$>$ In the case at least one of $\mathcal{R}_1$,\ $\mathcal{R}_2$ is of type\ I\hspace{-1.2pt}I\hspace{-1.2pt}I was studied in section 3.1\ and the case $\mathcal{R}_1$ is of type\ I was studied in \ref{10}.\ We studied the case at least one of $\mathcal{R} _1$,\ $\mathcal{R} _2$ is of type\ I\hspace{-1.2pt}I$_1$\ and the other is finite in \ref{28}.\ We studied the case either $\mathcal{R} _1$ is of 
type\ I\hspace{-1.2pt}I$_1$, $\mathcal{R} _2$ is not of type\ I\hspace{-1.2pt}I\hspace{-1.2pt}I and infinite,\ or 
$\mathcal{R} _1$ is of type I\hspace{-1.2pt}I$_\infty$ and $\mathcal{R} _2$ is of type\ $\mu$ ($\mu=$I, I\hspace{-1.2pt}I)
in \ref{37}.$\Box $

\vspace{\baselineskip}

\subsection{Other facts}
Next very beautiful proposition hard to come up with was proposed by Y.Ogata.
\begin{prop}\label{14}
Let $(\mathcal{R} ,\ \text{Ad}_{\Gamma }) $ be a graded\ von\ Neumann\
algebra. Then there is a projection\ $P \in Z(\mathcal{R})$ which satisfies following two conditions.\\
(i)\ For each projection\ $q \in Z(\mathcal{R})$ such that $q\leq P$ is even,\\
(ii) There is a projection\ $Q \in Z(\mathcal{R})$ such that \[\Gamma Q \Gamma =
(I-P)-Q. \]
\end{prop}
$<$proof$>$\ From \ref{13}, it suffices to see that proposition holds for $(Z(\mathcal{R}) ,\ \text{Ad}_{\Gamma})$. Thus we assume that $\mathcal{R}$ is commutative. Let 
$\{Q_a\}$ be a maximal orthogonal family of projections such that $\{Q_a+\Gamma Q_a \Gamma \} \subseteq \mathcal{R}^{(0)}$ is a orthogonal family of projections.\ Since $\{0\}$ satisfies the condition, the family is not empty. From Zorn's lemma, it has a maximal family. Let
\[P = I - \sum_{a} (Q_a + \Gamma Q_a \Gamma ). \]We shall show that this projection satisfies the conditions.\\
Let's see that $P$ satisfies (i).\ Let $q \leq P$. First, we consider the case $q \wedge \Gamma q \Gamma \neq 0$.\ Since 
$\mathcal{R}$ is commutative, 
$q \wedge \Gamma q \Gamma = q \Gamma q \Gamma $ and $q \Gamma q
\Gamma$ is even. Now,\ $q - q\Gamma q\Gamma $ is a projection which satisfies 
\[q -q\Gamma q\Gamma + \Gamma (q -q\Gamma q\Gamma)\Gamma = q + \Gamma
q\Gamma -2q \Gamma q \Gamma. \]\
We note that $q + \Gamma q\Gamma -2q \Gamma q \Gamma$ is even. Since $q + \Gamma q\Gamma -2q \Gamma q \Gamma$ self-adjoint and
\[(q + \Gamma q\Gamma -2q \Gamma q \Gamma)^2 = q + \Gamma q\Gamma -2q
\Gamma q \Gamma, \]
this is a projection.\ $P$ is also even because $\text{Ad}_\Gamma $ is strong-operator continuous. Thus it follows that 
$\Gamma q \Gamma \leq P$ for $q \leq P$ and $q + \Gamma q\Gamma -2q \Gamma q \Gamma
\leq q\vee \Gamma q \Gamma \leq P$. Consequently, if
$q \neq q \Gamma q \Gamma$, adjoining $q -q\Gamma
q\Gamma$ to $\{Q_a\}$ contradicts the maximality of $\{Q_a\}$.\ Thus $q= q \Gamma q \Gamma$ and $q$ is even.\\
Next we consider the case $q \wedge \Gamma q \Gamma = 0$.\ If $q \neq 0$, adjoining $q$ to $\{Q_a\}$ contradicts to the maximality of $\{Q_a\}$.\ If $q=0$, $q$ is even. We deduce that $q \leq P$ is even.\\
Next, we shall prove that $P$ satisfies (ii).\ Let $Q = \sum_{a} Q_a $ then
\[Q + \Gamma Q\Gamma = Q + \sum_{a} \Gamma Q_a \Gamma = I - P.
\]Thus $P$ satisfies (ii) .$\Box$
\vspace{\baselineskip}

\begin{prop}\label{41}
Let $(\mathcal{R} _1,\ \text{Ad}_{\Gamma _1})$ and $(\mathcal{R} _2,\ \text{Ad}_{\Gamma _2})$
be graded von Neumann algebras on $\mathcal{H} _1$,\ $\mathcal{H} _2$ respectively,\ and 
$A \in (\mathcal{R} _1)^{(0)}$,\ $B \in (\mathcal{R} _2)^{(0)}$.\ If the central support of 
$A\otimes B$ in $\mathcal{R} _1 \hat{\otimes } \mathcal{R} _2$
is $C$ and the central support of $A\otimes B$ in $\mathcal{R} _1\overline{\otimes
}\mathcal{R} _2$ is $D$, then $C=D$.
\end{prop}
$<$proof$>$ Let $X \in \mathcal{R} _1$,\ $Y \in \mathcal{R} _2$,\ $x \in
\mathcal{H} _1$ and $y \in \mathcal{H} _2$ then
\begin{align}
(X \otimes Y)(A \otimes B)(x \otimes y) = &(X \otimes Y^{(0)})(A
\otimes B)(x \otimes y)\\
 & + (X\Gamma _1 \otimes Y^{(1)})(A \otimes B)(\Gamma _1 x \otimes y)\\
 & \in C(\mathcal{H} _1\otimes \mathcal{H} _2).
\end{align}
Thus $D \leq C$. On the other hand,
\begin{align}
(X\Gamma _1 \otimes Y^{(1)})(A \otimes B)(x \otimes y) = &(X \otimes
Y^{(1)})(A \otimes B)(\Gamma _1 x \otimes y) \\
 & \in D(\mathcal{H} _1\otimes \mathcal{H} _2).
\end{align}
Thus $C \leq D$ and $C=D$.
$\Box $
\vspace{\baselineskip}

\begin{prop}\label{39}
Let $(\mathfrak{B} (\mathcal{H} ),\ \text{Ad}_{\Gamma })$ be a graded\ von\ Neumann\ algebra. Then there is a orthogonal family of minimal projections $\{E_a\}_a$ such that $E_a\in (\mathfrak{B}(\mathcal{H} ))^{(0)}$ and $\sum E_a=I$.
\end{prop}
$<$proof$>$ Let $\Gamma = P-(I-P)$ be spectral decomposition of $\Gamma $. We can choose CONS $\{x_a\}_a$ from $P(\mathcal{H} )$ and $(I-P)(\mathcal{H} )$. Let $E_a$ be a projection on the subspace generated by a single vector $x_a$. The family $\{E_a\}_a$ satisfies the condition.$\Box $
\vspace{\baselineskip}

\section{Acknowledgements}
This paper is English translation of master's thesis of the author. The author very thanks Professor Y.Ogata who is the advisor of the author for guiding this work.

\end{document}